\documentclass{article}
\usepackage{ijcai22}

\usepackage{times}
\usepackage{soul}
\usepackage{url}
\usepackage[hidelinks]{hyperref}
\usepackage[utf8]{inputenc}
\usepackage[small]{caption}
\usepackage{graphicx}
\usepackage{amsmath}
\usepackage{amsthm}
\usepackage{booktabs}
\usepackage{algorithm}
\usepackage{algorithmic}
\usepackage{algorithm}
\usepackage{algorithmic}
\usepackage{amssymb}

\usepackage{graphicx}
\usepackage{subfigure}
\usepackage{amsfonts}
\usepackage{amsmath,bm}
\usepackage{amssymb}
\usepackage{color}
\usepackage[round]{natbib}
\newtheorem{theorem}{Theorem}
\newtheorem{lemma}{Lemma}

\title{An Adaptive Incremental Gradient Method With Support for Non-Euclidean Norms}

\author{
Binghui Xie$^{*}$\footnote{$ :$ These authors contributed equally to this work }
\and
Chenhan Jin$^{*}$\and
Kaiwen Zhou\and
Wei Meng\and
James Cheng
\affiliations
The Chinese University of Hong Kong
}
\begin{document}

\maketitle

\begin{abstract}
Stochastic variance reduced methods have shown strong performance in solving finite-sum problems.
However, these methods usually require the users to manually tune the step-size, which is time-consuming or even infeasible for some large-scale optimization tasks.
To overcome the problem, we propose and analyze several novel adaptive variants of the popular SAGA algorithm.
Eventually, we design a variant of Barzilai-Borwein step-size which is tailored for the incremental gradient method to ensure memory efficiency and fast convergence. 
We establish its convergence guarantees under general settings that allow non-Euclidean norms in the definition of smoothness and the composite objectives, which cover a broad range of applications in machine learning.
The theoretical results supporting non-Euclidean norms fill the void of existing work.
Numerical experiments on standard datasets demonstrate a competitive performance of the proposed algorithm compared with existing variance-reduced methods and their adaptive variants.
\end{abstract}

\section{Introduction}
\label{s:section}
Many machine learning tasks involve solving the following optimization problem with a finite-sum structure:
\begin{equation}\label{problem1}
    \min_{x \in \mathbb{R}^d} \ f(x): = \frac{1}{n}\sum_{i \in [n]}f_i(x),
\end{equation}
where $x\in \mathbb{R}^d$ is the model parameter, each $f_i: \mathbb{R}^d \rightarrow \mathbb{R}$ is a loss function for the $i$-th training sample and $[n]:= \{1,\dots,n\}$.
More generally, the ``composite'' (or proximal) case covers a larger range of applications in machine
learning:
\begin{equation}\label{Problem}
    \min_{x \in \mathbb{R}^d} \ F(x) : = f(x) + h(x),
\end{equation}
where $h(\cdot): \mathbb{R}^d \rightarrow \mathbb{R}$ is a simple
and convex (but possibly non-differentiable) function, and the proximal operation of $h(\cdot)$ is easy to compute. %— few adaptive gradient methods are applicable in this setting.
Here, we also define $F_i(x) = f_i(x) + h(x)$ and $\nabla F_i(x) = \nabla f_i(x)+ \partial h(x)$ where $\partial h(x)$ denotes a sub-gradient of $h(\cdot)$ at $x.$
Throughout the paper, we denote $x^*$ as an optimal solution of Problem (\ref{Problem}).

However, when $n$ is very large, which is common in the modern learning tasks, gradient descent (GD) has a prohibitively high per-iteration cost.
%
%However, when nn is very large, GD has a prohibitively high per-iteration cost.
%
To mitigate the issue, stochastic gradient descent (SGD) \citep{10.1214/aoms/1177729586, DBLP:journals/siamjo/NemirovskiJLS09}, as an alternative, has been widely adopted to solve Problem (\ref{Problem}) in modern machine learning tasks.
Its simplest update process can be written as:
\begin{equation*}
    x_{k+1} = \arg \min_x \left\{\left\langle \widetilde{\nabla}_k, x \right\rangle  + \frac{1}{2\eta}\|x-x_k\|^2+ h(x)\right\},
\end{equation*}
where $\eta$ is the step-size, and the gradient estimator $\widetilde{\nabla}_k$ satisfies $\mathbb{E}[\widetilde{\nabla}_k] = \nabla f(x_k), \forall x \in \mathbb{R}^d$ to ensure the convergence.
The choice in SGD is to set $\widetilde{\nabla}_k=\nabla f_i(x)$.
However, due to the undiminishing variance $\mathbb{E}[\|\nabla f_i(x_k) - \nabla f(x_k) \|^2]$, SGD only converges at a sub-linear rate even if $F(\cdot)$ is strongly convex and smooth.

\paragraph{Variance reduction.}
Recently, some stochastic variance reduced (VR) methods have been proposed to solve the variance issue of SGD, such as SAG \citep{DBLP:conf/nips/RouxSB12}, SAGA \citep{DBLP:conf/nips/DefazioBL14}, SVRG \citep{DBLP:conf/nips/Johnson013}, SARAH \citep{DBLP:conf/icml/NguyenLST17} and SVRG-Loopless \citep{DBLP:conf/alt/KovalevHR20}.
%
% There are also several accelerated or generalized variants of these methods proposed these years, to list a few, \cite{DBLP:conf/nips/Defazio16,DBLP:journals/mp/GowerRB21,DBLP:journals/corr/abs-2006-11573,DBLP:conf/nips/ZhouSC20}.
%
The methods use better choices of the gradient estimator $\widetilde{\nabla}_k$ with its variance reducing as the algorithm converges.
Based on this property, in theory, VR methods only need $\mathcal{O}((n+\kappa)\log(\frac{1}{\epsilon}))$ stochastic gradient evaluations, while GD typically needs $\mathcal{O}(n\kappa \log(\frac{1}{\epsilon}))$ evaluations and SGD requires $\mathcal{O}(\frac{1}{\mu\epsilon})$ complexity.

Consequently, VR methods are commonly used in practice, especially for convex problems such as logistic regression \citep{DBLP:journals/corr/abs-2102-09645}.
However, all the VR methods above require a constant step-size, which depends on the characteristics of Problem (\ref{Problem}), such as the smoothness constant.
%
% In practice, these constants are often difficult to estimate.
%
In practice, researchers adopt a computationally expensive grid-search method to find a step-size.
%
%Besides such a computational cost, a 
Yet such constant step-sizes might lead to poor empirical performance.

\paragraph{Adaptive VR methods.}
Consequently, researchers have also proposed some adaptive VR methods to address the issue in recent years.
\citet{DBLP:conf/nips/TanMDQ16} proposed SVRG-BB based on Barzilai-Borwein (BB) step-size \citep{barzilai1988two} and derived SGD-BB and SAG-BB as heuristics.
\citet{c:21} developed SVRG-AS and SARAH-AS.
Although the methods mentioned above are capable of automatically changing step-size, they still introduce additional hard-to-tune hyperparameters.
\citet{DBLP:conf/icml/LiWG20} designed another variants of SVRG-BB and SARAH-BB given the knowledge of some problem-dependent constants.
More recently, \citet{DBLP:journals/corr/abs-2102-09645} proposed AdaSVRG for convex objectives, which extends the adaptive step-size in AdaGrad to SVRG's inner loop.
However, the theoretical analysis of AdaSVRG requires a bounded domain, and thus the algorithm needs to perform projection in each iteration, which is rarely the case in machine learning tasks.
\citet{DBLP:journals/corr/abs-2102-09700} emphasized the importance of adapting to the local geometry and proposed an implicit adaptive strategy for SARAH (AI-SARAH).
However, AI-SARAH requires solving a sub-problem in each iteration, which could be expensive for general objective functions other than linear regressions.

%In addition, most of the adaptive VR methods above are limited to solve Problem (\ref{problem1}) either in the convex setting or in the strongly convex setting.
%
%That is, AdaSVRG is only analyzed in the convex case, and other methods are analyzed in the strongly convex case.
%
%To the best of our knowledge, all the known adaptive VR methods and variants of SAGA do not work with a non-Euclidean norm in the definition of smoothness.
%
%However, allowing arbitrary norms contribute to wide applications in many areas of computer science \citep{allen2014linear}.

Although a considerable amount of work has been done for adaptive SVRG and SARAH, another stochastic variance reduction method, SAGA, does not have any adaptive variant. 
SAGA, which is free of inner loop length tuning and usually has a better practical performance compared with SVRG, has been implemented in the scikit-learn package \citep{scikit-learn} (while SVRG and SARAH have not been) to solve linear regression tasks.
Therefore, an adaptive variant of SAGA is of great interest.

However, unlike SVRG and SARAH, which conduct a full gradient computation in each outer iteration and match the adaptive step-size frameworks like BB step-size \citep{DBLP:conf/nips/TanMDQ16},
SAGA does not compute full gradients. Instead, SAGA maintains a table of the component function gradients $\nabla f_i(\cdot)$ that were lastly evaluated. Naively using the strategies of the adaptive SVRG or SARAH will not make full use of incremental gradients, hence will introduce a high computation cost for SAGA, resulting in a poor memory overhead. In addition, there exists a discrepancy in the theory between VR methods and their adaptive variants. VR methods are proven to be able to solve Problem \ref{problem1} and \ref{Problem} both in convex and strongly convex settings. However, most of the adaptive VR methods are limited to solving Problem \ref{problem1} either in the convex setting or in the strongly convex setting.
%For instance, AdaSVRG is only analyzed in the convex case, and other methods are analyzed in the strongly convex case.
Moreover, although Euclidean norm space is common and widely considered in both academia and industry, there exist some problems related to non-Euclidean norms that are also crucial to study both in theory and industry, see Appendix \ref{appendix: importance} for details. %Some VR methods can be modified to work with non-Euclidean norms in the definition of smoothness, e.g., Kayusha in \citep{allen2014linear}. 
To the best of our knowledge, all the known adaptive VR methods and variants of SAGA do not work with a non-Euclidean norm. Such challenges lead us to the natural question: \textit{is it possible to design an efficient adaptive variant for SAGA that can address the above problems under theoretical guarantees?}
%

%Unlike SVRG and SARAH, {the adaptive variants have been properly addressed in [], it is unknown to have a adaptive variant for SAGA. The strategies of BB,  }which conduct a full gradient computation in each outer iteration and match the adaptive step-size framework like BB step-size \citep{DBLP:conf/nips/TanMDQ16},
%SAGA does not compute full gradients.
%
%Instead, SAGA maintains a table of the component function gradients $\nabla f_i(\cdot)$ that were lastly evaluated.
%
% Then SAGA can get rid of the full gradient computation.
%
%Naively using the strategies of the adaptive variants of SVRG and SARAH will not make full use of incremental gradients, hence will introduce a higher computation overhead for SAGA. 
%
%Therefore, it remains as a challenge to design an efficient adaptive variant of SAGA.
%{\color{red} Such a challenge leads us to the natural question \textit{is it possible to design an efficient adaptive variant for SAGA under statistical guarantees? }} 
\subsection{Contribution}

In this paper, we give a positive answer to above question. We propose the first adaptive framework of SAGA, which we call SAGA-BB, to solve Problem \ref{problem1} and \ref{Problem} both in (strongly) convex and non-Euclidean norm space settings .
%
% Moreover, we do not assume the smoothness of the objective is defined on $\ell_2$-norm, which requires a new and more challenging analysis.
%
Specifically, we summarize our contributions as follows.
\begin{itemize}
    \item We propose a novel adaptive stochastic incremental gradient framework based on SAGA.
    With this framework, we design some novel variants of the existing adaptive step-size, which make full use of the incremental gradients maintained by SAGA.
    These variants can automatically adjust to the local geometry.
    \item We also propose a variant of BB step-size.
    This leads to SAGA-BB, a simple yet powerful adaptive VR method.
    SAGA-BB is almost tune-free and presents more robust and consistent practical performance than its counterparts.
    Moreover, SAGA-BB enjoys the same efficient implementation for linear models as the original SAGA method, while allowing both single sample and mini-batch settings.
    \item We conduct convergence analysis on the composite problem which supports a proximal operator.
    We prove that SAGA-BB achieves a linear convergence rate $\mathcal{O}((n+\kappa)\log(\frac{1}{\epsilon}))$ for strongly convex problems, and a sub-linear rate $\mathcal{O}(\frac{1}{\epsilon})$ for general convex objectives.
    Both rates match the oracle complexity achieved by SAGA \citep{DBLP:conf/nips/DefazioBL14}.
    Our theoretical results allow non-Euclidean norms in the definition of smoothness.
    This fills the void of the existing analyses of SAGA.
    \item We conduct experiments on common machine learning tasks to verify the robustness and effectiveness of SAGA-BB.
    SAGA-BB demonstrates a strongly competitive performance compared with fine-tuned SAGA, SVRG, SARAH, SVRG-BB and other related methods, across different standard datasets.
    %\WM{Can we say it outperforms the others?}
    %\BH{Sure, SAGA-BB indeed outperforms the others in our experiments. I was afraid some reviews were the authors of these papers, so I used a mild statement. I have used the "superior performance" to replace the original one}
\end{itemize}

\section{Algorithm}
\label{s:algorithm}
In this section, we introduce our adaptive variant of SAGA, named SAGA-BB. 
%and establish its memory and performance efficiency.
%
Before presenting our results, we first discuss the existing adaptive step-sizes in the deterministic setting (Section \ref{s:algorithm-deterministic}), which inspire our work.
We next propose several stochastic adaptive step-sizes for SAGA and demonstrate SAGA-BB is the best among them (Section \ref{s:algorithm-saga-bb}).

\subsection{Deterministic Adaptive Step-size}
\label{s:algorithm-deterministic}
We summarize several successful adaptive step-sizes for gradient methods in the deterministic setting, which inspired us.

\textbf{Barzilai-Borwein step-size.}
Barzilai-Borwein step-size is inspired by quasi-Newton methods, which solves Problem (\ref{problem1}) by choosing $\eta$ as an approximation of $H^{-1}_t$ for GD, where $H$ is the Hessian matrix of $f$.
However, one needs to solve a linear system to obtain such an approximation, which is time consuming in some cases.
To solve the problem, \citet{barzilai1988two} proposed BB step-size to satisfy the residual of the secant equation, which can be formulated as:
$\min\  \|\eta_k^{-1}s_k-y_k \|$ or $ \min\ \|s_k-\eta_k y_k  \|,$
where $s_k = x_k - x_{k-1}$ and $y_k = \nabla f(x_k) - \nabla f(x_{k-1}).$ The corresponding solutions are
\begin{equation*} \label{BBmethods}
    \eta_k^{BB1} = \frac{s_k^Ts_k}{s_k^Ty_k} \ \ \ \ or \ \ \ \ \eta_k^{BB2} =  \frac{s_k^Ty_k}{y_k^Ty_k},
\end{equation*}
respectively.
%
% Such two schemes perform similarly for convex quadratic functions \citep{barzilai1988two}.
%
% Most papers afterwards prefer to mention and analyze only the first scheme, BB1.
%
% However, in our experiments, we notice that BB2 can have a better and more stable performance in the stochastic setting (Section \ref{s:algorithm-saga-bb}).

\textbf{ALS step-size.}
It is well known that with the constant step-size $\eta \leq 1/L$, GD has a sub-linear rate for convex and smoothness objectives, where $L$ is the smoothness constant.
To obtain an estimation of the local smoothness constant, \citet{pjm/1102995080,vrahatis2000class, c:21} proposed the ALS step-size. That is,
\begin{equation*}
    \eta_k^{ALS} = \frac{\|x_k -x_{k-1} \|}{2\|\nabla f(x_k)- \nabla f(x_{k-1}) \|}.
\end{equation*}
With the scheme above, ALS step-size can obtain the local estimation of the smoothness constant $1/L$ adaptively without additional information.

\textbf{YK step-size.}
\citet{DBLP:conf/icml/MalitskyM20} proposed another $L$ estimation method---YK step-size:
\begin{equation*}
    \eta_k^{YK} = \min \left \{\sqrt{1+\theta_{k-1}}\eta_{k-1}, \frac{\|x_k -x_{k-1} \|}{2\|\nabla f(x_k)- \nabla f(x_{k-1}) \|} \right \}.
\end{equation*}
The $\theta_k$ in YK step-size and its variants introduced below (Section \ref{s:algorithm-saga-bb}) is defined as $\eta_k/\eta_{k-1}$.
This adaptive step-size leads to a novel analysis technique for gradient descent.
From the perspective of the scheme, we can regard the YK step-size as a bounded variant of ALS step-size.

% \textbf{Polyak’s step-size.}
% \citet{polyak1987introduction} proposed a step-size that largely simplifies the convergence analysis. The Ployak's step-size at iteration $k$ is: $\frac{f(x_k)-f^*}{\|\nabla f(x_k) \|^2}.$
% Polyak’s step-size stands out in the deterministic settings.
% %
% \citet{DBLP:conf/aistats/LoizouVLL21} proposed a stochastic variant of Polyak's step-size, which requires the knowledge of $f_i^*$, which is hard to know in some cases.
% %
% Thus, we do not consider this method in our work.

Note that all the above methods require a full gradient evaluation per epoch or per iteration.
Directly extending them to the incremental gradient method is impractical.
To tackle these issues, we next provide several novel variants to incorporate these methods to SAGA.
%\TODO{Add a cross reference to the later section, unless your analysis follows this sub-section immediately.}
%\BH{The next sub-section is going to give the analysis.}

\subsection{Adaptive Step-sizes for SAGA}
\label{s:algorithm-saga-bb}
%\WM{From the early introduction, I was expecting that you would discuss some other methods in this setting.
%Is your method the only stochastic adaptive step-size scheme? If yes, you might want to change the sub-section title and rephrase in the introduction paragraph.}
%\BH{Thank you. I have modified the early introduction. (Does it denote the first paragraph in this section?) Due to the limited pages, and no such methods for SAGA, we just discuss our methods.}
%
We propose a tailor-made adaptive step-size scheme for SAGA, presented in Algorithm \ref{alg:stc schm}.
Since we want to adjust the step-sizes based on the local geometry of our objective, it is natural to capture the change in the points and gradients with $s_k$ and $y_k$ respectively, in each iteration.
With the local information, the step-size will be adjusted automatically every $m$ iterations.
Since there are $n$ points, we set the frequency $m$ as $n$ in practice.
Although the scheme requires an initial step-size, we observe from numerical experiments that the performance of our method is not sensitive to the choice of the initial step-size.

\begin{algorithm}[tb]
\caption{SAGA with adaptive step-size scheme}
\label{alg:stc schm}
\textbf{Input}: Max epochs $K$, initial point $x_1$, initial step size $\eta_1$, update frequency $m$.  \\
\textbf{Initialize}: Initialize the gradient table $\nabla f_i(\phi_i^{1})$ with $\phi_i^{1}=x_{1}$, and the average $\mu_1 = \frac{1}{n}\sum_{i=1}^n \nabla f_i(\phi_i^{1})$.  \\
\begin{algorithmic}[1] %[1] enables line numbers
\FOR{$k = 1,2,\dots, K$}
\IF {$k \mod m == 0$}
\STATE
update $\eta_k$
\ELSE
\STATE
$\eta_k = \eta_{k-1}$
\ENDIF

\STATE Randomly choose a sample $i$ and take $\phi_i^{k+1}=x_{k}$, store $\nabla f_{i}(\phi_{i}^{k+1})$ in the table and keep other entries unchanged.
\STATE
$\widetilde{\nabla}_{k} =\nabla f(\phi^{k+1}_i)-\nabla f(\phi^{k}_i)+\mu_{k}$
\STATE
$y_{k}= \nabla f_{i}(\phi_{i}^{k+1})-\nabla f_{i}(\phi_{i}^{k})$
\STATE
$s_{k}= \phi_{i}^{k+1}-\phi_{i}^{k}$
\STATE
$\mu_{k+1}=\mu_{k} + \frac{1}{n}(\phi^{k+1}_i-\phi^{k}_i)$
\STATE 
$x_{k+1}=\arg \min_x \Big\{\left\langle \widetilde{\nabla}_{k}, x \right\rangle  + \frac{1}{2\eta_k}\|x-x_k\|^2+ h(x)\Big\}$

%$x_{k+1} = x_{k} - \eta_k \widetilde{\nabla}_{k}$

\ENDFOR
\STATE \textbf{return} $x_{K+1}$
\end{algorithmic}
\end{algorithm} 

Based on our framework (Algorithm \ref{alg:stc schm}), it is natural to update the step-size by the arithmetic mean of the information from the last $m$ iterations.
We propose the following variants of the deterministic adaptive step-sizes in Section \ref{s:algorithm-deterministic} for SAGA.

\textbf{The Stochastic Variant of BB1 step-size.}
\begin{equation*}
    \frac{1}{m}\sum_{i=k-m}^{k}\frac{1}{\alpha}\frac{s_i^Ts_i}{s_i^Ty_i}.
\end{equation*}

\textbf{The Stochastic Variant of BB2 step-size.}
\begin{equation*}
    \frac{1}{m}\sum_{i=k-m}^{k}\frac{1}{\alpha}\frac{s_i^Ty_i}{y_i^Ty_i}.
\end{equation*}

\textbf{The Stochastic Variant of ALS step-size.}
\begin{equation*}
    \frac{1}{m}\sum_{i=k-m}^{k}\frac{\|s_i \|}{\alpha\|y_i\|}.
\end{equation*}

\textbf{The Stochastic Variant of YK step-size.}
\begin{equation*}
    \min \left  \{\sqrt{1+\theta_k} \eta_{k-1},  \frac{1}{m}\sum_{i=k-m}^{k}\frac{\|s_i \|}{\alpha\|y_i\|}\right \}.
\end{equation*}

The constant $\alpha$ in the above stochastic variants is a common hyperparameter to control the convergence, which is similarly required in previous works \citep{DBLP:conf/icml/LiWG20,DBLP:conf/aistats/LoizouVLL21,DBLP:conf/nips/TanMDQ16}. 

\begin{figure}[t]
    \centering
    \subfigure[]{
    \includegraphics[width=0.20\textwidth]{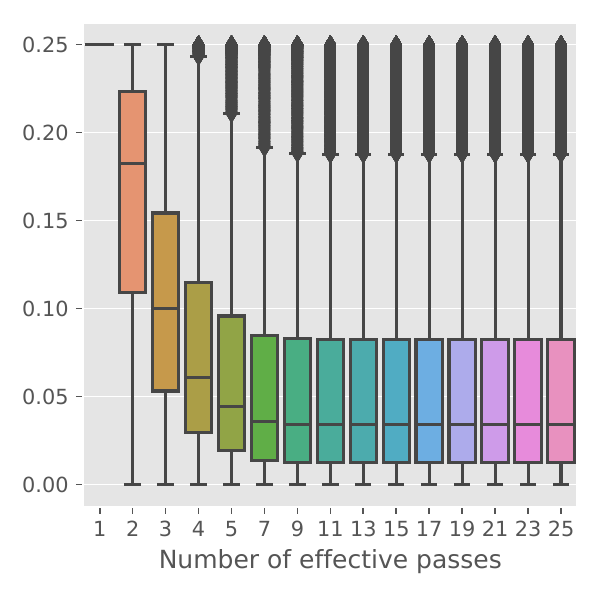}
    \label{fig: local information-1}}
    \subfigure[]{
    \includegraphics[width=0.20\textwidth]{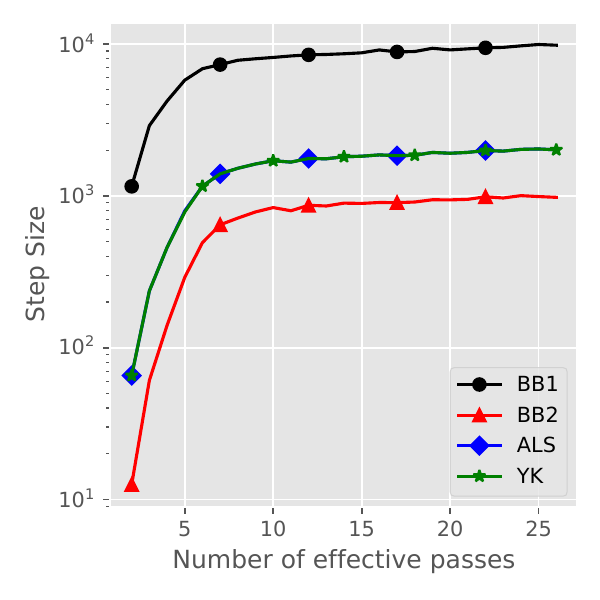}
    \label{fig: local information-2}}
    \caption{The distribution of $p_i$ on the left and the evolution of the stochastic step-sizes on the right.}
    \label{fig: local information}
\end{figure}

To demonstrate how the above adaptive step-sizes are related to local information, we conduct a heuristic experiment on the \emph{ijcnn} dataset.
% $\footnote{\mbox{The dataset is available at https://www.csie.ntu.} edu.tw/˜cjlin/libsvmtools/datasets/.}$
%
Inspired by \cite{DBLP:journals/corr/abs-2102-09700}, we consider Problem (\ref{problem1}) with each $f_i(x)=\log(1+\exp(-b_ia_i^Tx)+\lambda\|x\|^2$ as an example, where $a_i \in \mathbb{R}^d$ is the data vector, and $b_i \in \{-1, +1\}$ is the label.
The local curvature of $f(x)$ can be denoted as its Hessian:$\nabla^2f(x) = \frac{1}{n}\sum_{i=1}^{n}\frac{\exp(-b_ia_i^Tx)}{[1+\exp(-b_ia_i^Tx)]^2}a_ia_i^T+\lambda I.$ It is easy to derive that the local curvature is highly related to the $p_i(x)=\frac{\exp(-b_ia_i^Tx)}{[1+\exp(-b_ia_i^Tx)]^2}$.
We plot the distribution of $p_i$ in Figure \autoref{fig: local information-1}.
Correspondingly, we plot the evolution of step-sizes generated by these variants in Figure \autoref{fig: local information-2}.
The constant $\alpha$ is set to 1 here.
%
% We could see that $p_i$ is decreasing as the number of passes increases, and stabilizes after around 20 passes while all the step-sizes increase as the number of passes grows.
We could see that $p_i$ is decreasing as the number of passes increases, and stabilizes after around 20 passes.
Thus, the local curvature will be smaller as well and a larger step-size can be used.
Correspondingly, it is reasonable that the step-sizes increase as the number of passes grows, as shown in Figure \autoref{fig: local information-2}. 

However, the trajectories of these step-sizes are much steeper than that of $p_i$, which can lead to unstable performance of algorithms and divergence.
The reason is that arithmetic mean scheme is highly affected by the outliers and sensitive to the difference of $f_i$.
Take the variant of BB1 step-size as an example. $s_i^Ty_i$ can be much smaller for some indexes, which is also found in \cite{fletcher2005barzilai,burdakov2019stabilized}.
This causes the $\frac{s_i^Ts_i}{s_i^Ty_i}$ to have a very large value, which consequently leads to a large step-size.

To reduce the impact of the ``outliers'', we use the mediant instead of arithmetic mean to gather the information.
The mediant can be regarded as weighted arithmetic mean, where the weight lies in the denominator.
Therefore, for the ``outliers'' with small denominators, their weights would also be small.
They will contribute less to the update of step-size.
Based on these considerations, we propose the stable adaptive step-size for SAGA as follows:

\textbf{Stable BB1 step-size.}
\begin{equation*}
    \frac{\sum_{i=k-m}^{k} s_{i}^Ts_{i}}{\alpha \sum_{i=k-m}^{k} s_{i}^Ty_i}.
\end{equation*}

\textbf{Stable BB2 step-size.}
\begin{equation*}
    \frac{\sum_{i=k-m}^{k} s_{i}^Ty_{i}}{\alpha \sum_{i=k-m}^{k} y_{i}^Ty_i}.
\end{equation*}

\textbf{Stable ALS step-size.}
\begin{equation*}
    \frac{\sum_{i=k-m}^{k} \|s_{i}\|}{\alpha \sum_{i=k-m}^{k} \|y_i\|}.
\end{equation*}

\textbf{Stable YK step-size.}
\begin{equation*}
    \min \left  \{\sqrt{1+\theta_{k-1}}\eta_{k-1},  \frac{\sum_{i=k-m}^{k} \|s_{i}\|}{\alpha \sum_{i=k-m}^{k} \|y_i\|}\right \}.
\end{equation*}

\begin{figure}[t]
    \centering
    \subfigure[]{
    \includegraphics[width=0.20\textwidth]{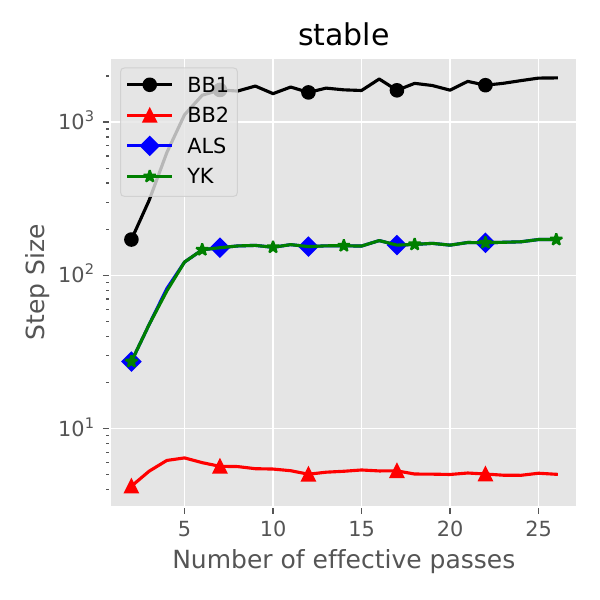}
    \label{fig: stable methods-1}}
    \subfigure[]{
    \includegraphics[width=0.20\textwidth]{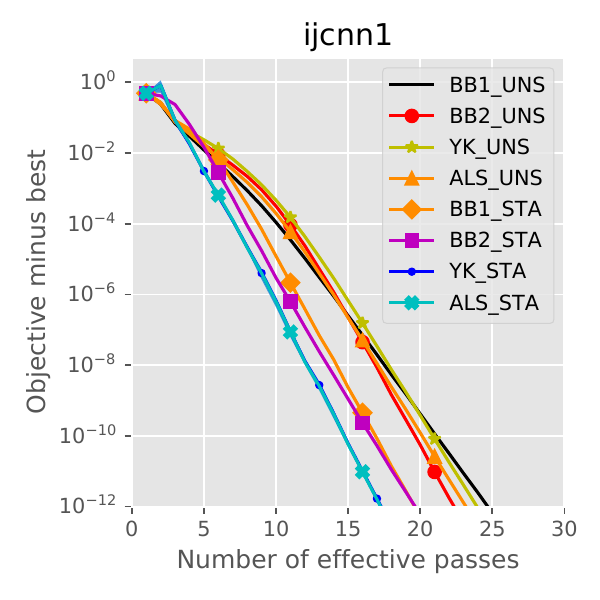}
    \label{fig: stable methods-2}}
    \caption{The evolution of step-sizes on the left and the performance for different adaptive methods.}
    \label{fig: stable methods}
\end{figure}

We now refer the previous stochastic adaptive step-sizes as the unstable ones. 
We plot the evolution of these stable variants in Figure \autoref{fig: stable methods-1}.
From the result, we can see the stable step-sizes show a much smaller and smoother performance than those unstable ones in Figure \autoref{fig: local information-2}.
For illustration, we also compare the performance of the stable methods with the unstable methods on the \emph{ijcnn} dataset.
Since $\frac{\|s_i \|}{\|y_i\|}$ is the geometrical mean of $\frac{s_i^Ts_i}{s_i^Ty_i}$ and $\frac{s_i^Ty_i}{y_i^Ty_i}$, we set $\alpha=2$ for stable BB2 step-size, $\alpha=\sqrt{m/2}$ for stable YK and ALS methods, and $\alpha=m$ for stable BB1 step-size across different datasets and losses.
For the \emph{ijcnn} dataset, we tuned $\alpha=20/m$, $\alpha=5/m$, $\alpha=5/m$, $\alpha=4/m$ for unstable BB2, YK, ALS, BB1 step-sizes, respectively.
Note that different from unstable step-sizes, these choices for stable methods also work well on other datasets (refer to Appendix \ref{appendix: exps}).
As shown in Figure \autoref{fig: stable methods-2}, the stable methods all enjoy a faster convergence.

Among these stable methods, BB2 shows a more consistent performance across different datasets (refer to the Appendix \ref{appendix: exps}).
In addition, as implemented in Appendix \ref{appendix: implement}, BB2 step-size only requires $\mathcal{O}(n)$ memory consumption for problems with a linear predictor , while other methods need $\mathcal{O}(nd)$ memory overhead.
Therefore, the stable BB2 step-size is the sensible choice for SAGA.

\section{Theory in General Norms}
In this section, we theoretically analyze the performance of SAGA-BB in general norm spaces and composite cases.
We first introduce some preliminaries about general norms.
\subsection{Preliminary}
\textbf{Norm Space and Notations.}
Let $\mathbf{E}$ be a real vector space with finite dimension $d$ and $\mathbf{E}^*$ is its dual space.
$\langle y,x \rangle$ represents the value of a linear function $y \in \mathbf{E}^*$ at $x \in \mathbf{E}.$
$\|\cdot\|$ denotes an arbitrary norm over $\mathbf{E}.$
The dual norm $\|\cdot\|_*$ over the dual space is defined in the standard way:
%\begin{equation}
    $\|y\|_*:= \max\{\langle y,x \rangle:\ \|x\| \leq 1\}.$
%\end{equation}
For instance, $\ell_p$ norm is the dual of $\ell_q$ norm if $\frac{1}{p}+\frac{1}{q}=1.$
\\
\\
\textbf{Non-Euclidean Norm Smoothness.} We consider the smoothness with respect to an arbitrary norm $\|\cdot\|$ over $\mathbf{E}.$
Symbolically, we say the differentiable function $f(\cdot)$ is $L$-smooth w.r.t. $\|\cdot\|$ if $\forall x,y \in \mathbf{E},$ it satisfies: $\|\nabla f(x)-\nabla f(y)\|_* \leq L\|x-y\|.$ Additionally, if $f(\cdot)$ is convex, the following holds \citep{DBLP:books/sp/Nesterov04}: $\|\nabla f(x)-\nabla f(y)\|_* \leq 2L(f(x) -f(y) -\langle \nabla f(y),x-y \rangle).$ Some famous problems have better smoothness parameters when non-Euclidean norms are adopted. See the discussions in \cite{allen2014linear}. Moreover, We point out in Appendix \ref{appendix: importance} that it is important for SAGA to support non-Euclidean norms.

\subsection*{Bregman Divergence}
%\WM{Add an attributive to the setting.}
We introduce the proximal setting, which generalizes the usual Euclidean setting \citep{zhou2020amortized}.
%Following the traditions in the non-Euclidean norm setting \cite{allen2014linear}, we select
%We introduce \XXX{the setting} \cite{allen2014linear} which generalizes the usual Euclidean setting.
The distance generating function (DGF) $d: \mathbf{E} \rightarrow \mathbb{R}$ is required to be continuously differentiable and 1-strongly convex with respect to $\|\cdot\|,$ i.e., $d(x)-d(y) - \langle \nabla d(y),x-y \rangle \geq \frac{1}{2}\|x-y\|^2,\ \forall x,y \in \mathbf{E}.$
%\WM{Is the above an inline equation?}
%\BH{Yes, we have revise this.}
%
Accordingly, the prox-term (Bregman divergence) is given as
\begin{equation*}
    V_d(x,y) \overset{\text{def}}{=} d(x)-d(y) - \langle \nabla d(y),x-y \rangle,\ \forall x,y \in \mathbf{E}.
\end{equation*}
The property of $V_d(\cdot,\cdot)$ ensures that $V_d(x,x)=0$ and $V_d(x,y) \geq
\frac{1}{2}\|x-y\|^2 \geq 0,\ \forall x,y \in \mathbf{E}.$
One of the key benefits that comes from such a setting is: by adjusting $\|\cdot\|$ and $d(\cdot)$ to the geometry of the problem, mirror descent achieves a smaller problem-dependent constant than the Euclidean algorithms \citep{nemirovskij1983problem}.

We first present two instances of Bregman Divergences:
%\WM{Try to rephrase the following two items. They read grammatically incorrect to me.}
%\BH{Yes, there are some grammatical errors. We have revised them.}

\begin{itemize}
    \item With $\|\cdot\|$ being $\|\cdot\|_2,$ $\mathbf{E}=\mathbb{R}^d,$  $d(x)=\frac{1}{2}\|x\|^2_2,$ this setting can be taken as the standard Euclidean setting: $V_d(x,y)=\frac{1}{2}\|x-y\|^2_2.$
    \
    \item With $\|\cdot\|$ being $\|\cdot\|_1,$ $d(x)=\sum_i x_i\log x_i,$ this setting can be taken over simplex $\mathbf{E} \subseteq \Delta^d \overset{\text{def}}{=}\{x \in \mathbb{R}_+^d:\sum_{i=1}^d x_i =1\},$ known as the entropy function:
    %\begin{equation}
      $ V_d(x,y)=\sum_i x_i\log x_i/y_i \geq \frac{1}{2}\|x-y\|_1^2.$
    %\end{equation}
\end{itemize}

We select the $V_d$ such that the proximal operator, $\text{Prox}_h(x,\mathcal{G}) \overset{\text{def}}{=}\arg \min_{u \in \mathbf{E}} \left\{V_d(u,x)+\left\langle \mathcal{G} , u \right\rangle+h(u)\right\}$ is easy to compute for any $x \in \mathbf{E}, \mathcal{G} \in \mathbf{E}^*.$ Examples can be found in \cite{parikh2014proximal}.

\textbf{Algorithm Changes.} Suppose $f(x)$ is $L$-smooth w.r.t. $\|\cdot\|$ and a bregman divergence $V_d(x,y)$ is given, we perform the following change to the algorithms. In line 12 of Algorithm \ref{alg:stc schm}, we change the $\arg \min$ (proximal operation) to be its non-Euclidean norm variant \citep{DBLP:conf/stoc/Zhu17}: $x_{k+1}= \text{Prox}_h(x_k,\eta_k \widetilde{\nabla}_k) .$

\textbf{Generalized Strong Convexity of $\bm{f}$.} A convex function $f$ is strongly convex with respect to $ V_d(x,y)$ rather than the $\|\cdot\|,$ i.e., $f(x) \geq f(y) + \langle \nabla f(y),x-y \rangle +\mu V_d(x,y) ,\ \forall x,y \in \mathbf{E}.$

This can be satisfied if $d(\cdot) \overset{\text{def}}{=}\frac{1}{\mu}f(\cdot).$
This is known as the ``generalized strong convexity'' \citep{shalev2007online} and is necessary for the existing linear convergence results in the strongly convex setting.
Note that in Euclidean setting, ``generalized strong convexity'' is equivalent to strong convexity that is with respect to $\ell_2$ norm.

\subsection{Convergence Analysis}
In this section, we analyze the convergence of SAGA-BB.
%All expectations are taken with respect to the choices of the sample ii at current iteration kk and conditioned on xkx_k and ∇fi(ϕki)\nabla f_i(\phi_i^k) unless otherwise specified.
%
Note that the analysis below can be extended to the mini-batch version with a random sample size $|B|$ (see Appendix \ref{appendix: mini-batch} for details).
First, we present the important lemma which provides a variance bound of the stochastic gradient estimator used by SAGA-BB.

\begin{lemma} \label{VB}
	(Variance Bound). For the stochastic gradient estimator $\widetilde{\nabla}_k$ in Algorithm \ref{alg:stc schm}, we have the following variance bound:
	\begin{align*}
			&\mathbb{E}\left[\left\|\nabla f(x_k)-\widetilde{\nabla}_{k}\right\|^2_*\right]\\
			\leq&{} 4L\left[f(x_k)-f(x^*)- \langle \nabla f(x^*), x_k - x^*\rangle \right]\\
			+&4L\left[\frac{1}{n}\sum_{i=1}^{n}f_i(\phi_i^k)-f(x^*)- \frac{1}{n}\sum_{i=1}^{n}\langle \nabla f_i(x^*), \phi_i^k - x^*\rangle \right].
	\end{align*}
\end{lemma}

Now we can formally present the main theorem of SAGA-BB below.
Compared with other analyses of incremental methods  \citep{DBLP:conf/nips/DefazioBL14,zhou2019direct} that were built only in Euclidean space, we show that SAGA-BB can operate in non-Euclidean norms. Explicitly,  we design a new Lyapunov function, which differs from \cite{DBLP:conf/nips/DefazioBL14} and helps us extend our theorems into general norms. For the general convex case, we establish the convergence result of SAGA-BB in terms of averaged iterate: $\bar{x}_K = \frac{1}{k}\sum_{i=1}^k x_i.$

\begin{theorem} \label{theo1}
	(General Convex). Define the Lyapunov function $T,$ which modifies the one in \cite{DBLP:conf/nips/DefazioBL14}:
	\begin{align*}
	    &T_k:=\\
	    &\frac{1}{\omega}\left[\frac{1}{n}\sum_{i=1}^{n}F_i(\phi_i^k)- F(x^*)-\frac{1}{n}\sum_{i=1}^{n}\langle \nabla F_i(x^*),\phi_i^k - x^* \rangle \right].
	\end{align*}
	Then we have that for $\eta_k \leq \frac{1}{9L}$ and a constant  $\alpha > 9$:
	\begin{align*}
    \mathbb{E}[F(\bar{x}_{K+1})] - F(x^*) \leq& \frac{2}{K}\bigg[\frac{n}{4}(F(x_1)-F(x^*))\\
    &{}+\frac{1}{2}F(x_1)+\alpha L V_d(x^*,x_1)\bigg].
    \end{align*}
\end{theorem}

Thus, for general convex objectives, SAGA-BB yields a sub-linear convergence rate $\mathcal{O}( \frac{1}{\epsilon}),$ which keeps up with the oracle complexity achieved by SAGA \citep{DBLP:conf/nips/DefazioBL14}.

Moreover, unlike SAGA, we do not need to assume that Problem (\ref{Problem}) is considered on Euclidean norm space, which makes our analysis more general.
Compared with other adaptive methods, SAGA-BB can support composite problems where the proximal operator is used.

\subsection{Strongly Convex Case}
We consider $f$ satisfies the generalized strongly convexity. For this case, we modify Algorithm \ref{alg:stc schm} to a new method called R-SAGA-BB. In each iteration, R-SAGA-BB runs Algorithm \ref{alg:stc schm} for $K$ iterations from the current point $\bar{x}_K$ and uses its output $\bar{x}_{K+1}$ as the starting point for the next iteration. This strategy is known as the restarting technique \citep{dvurechensky2016stochastic,gorbunov2020stochastic}. By choosing $K$ properly, we get an accelerated method for the generalized strongly convexity of $f.$

\begin{theorem} \label{theo2}
    (Strongly Convex). Assuming $f$ is generalized strongly convex, if we choose $K = \frac{\mu n + 6\alpha L}{\mu}$ with a constant $\alpha > 9$, then we have after $\tau$ times running Algorithm \ref{alg:stc schm}, 
    \begin{equation*}
        \mathbb{E}[F(\bar{x}_{\tau+1})] - F(x^*) \leq \frac{1}{2^{\tau}}(F(x_1) - F(x^*)).
    \end{equation*}
\end{theorem}

Therefore, for the generalized strongly convex objective, although we adopt a restarting technique, SAGA-BB yields a fast linear convergence rate $\mathcal{O}((n+\frac{L}{\mu})\log \frac{1}{\epsilon})$ which keeps up with the oracle complexity achieved by SAGA.

\begin{figure*}[ht]
    \centering
    \subfigure[]{
    \includegraphics[width=0.20\textwidth]{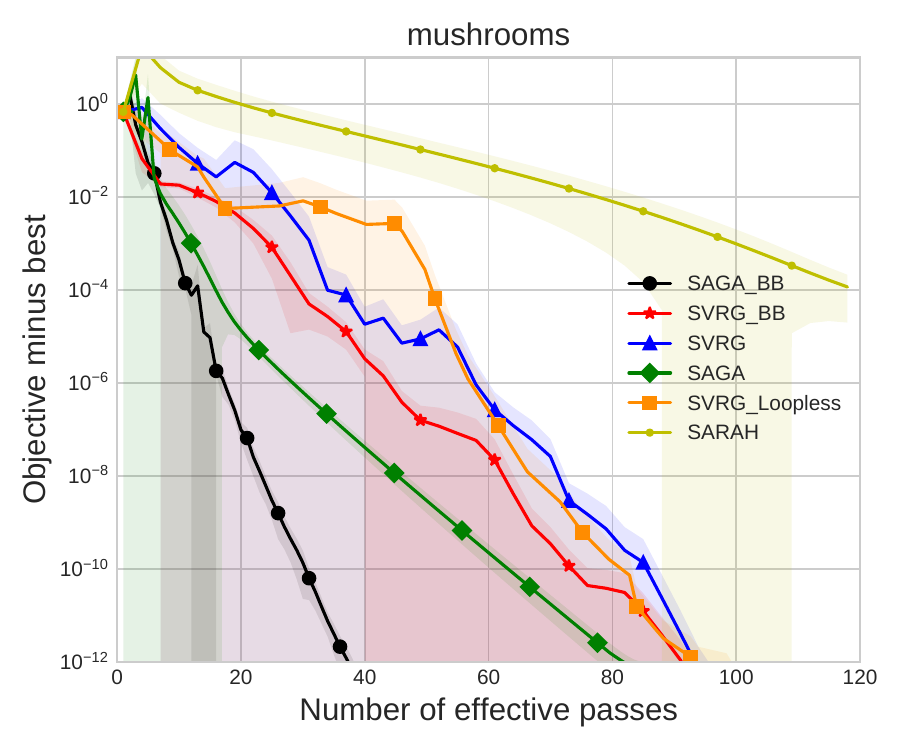}}
    \subfigure[]{
    \includegraphics[width=0.20\textwidth]{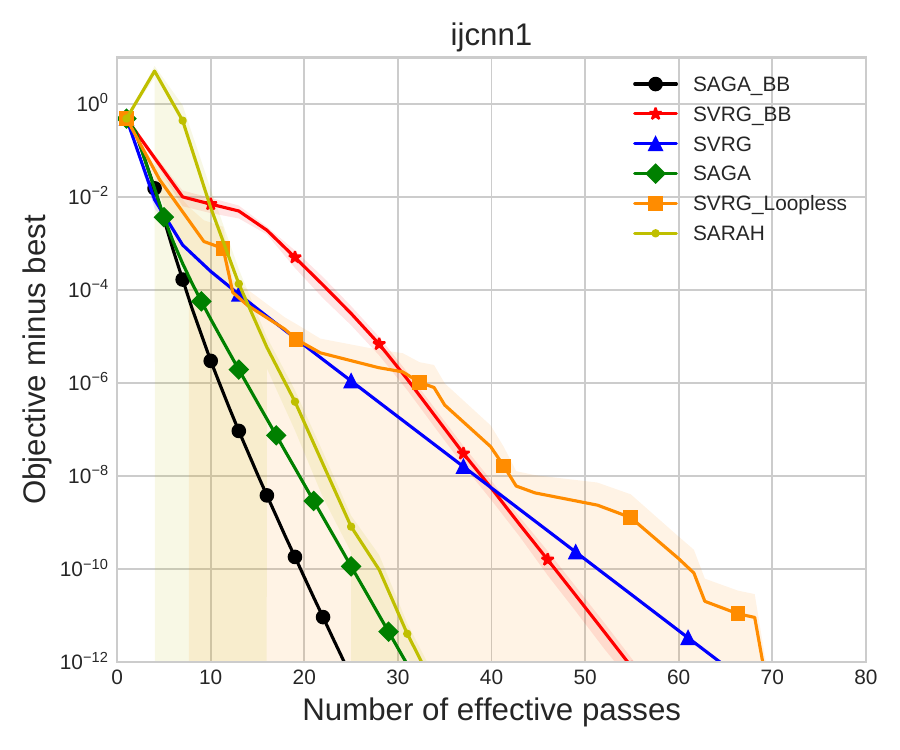}}
    \subfigure[]{
    \includegraphics[width=0.20\textwidth]{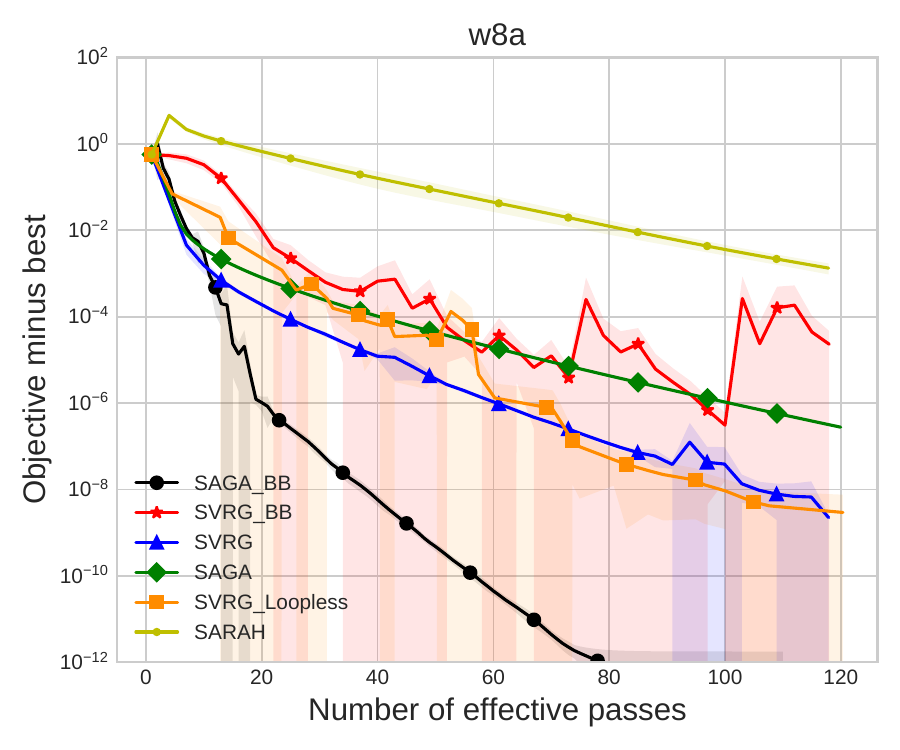}}
    \subfigure[]{
    \includegraphics[width=0.20\textwidth]{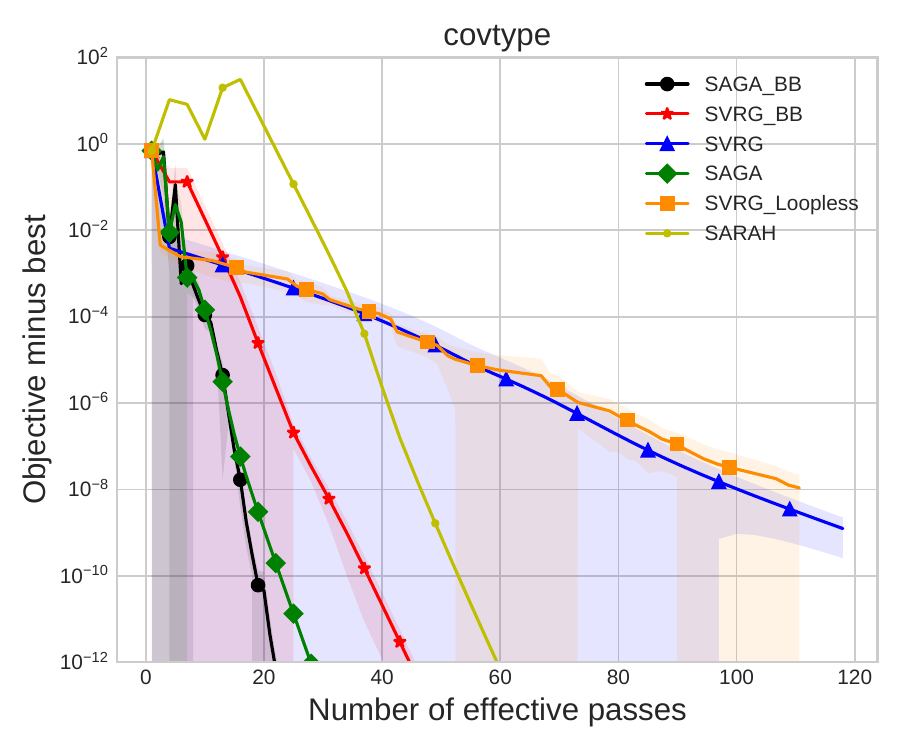}}
    \caption{Convergence on different LibSVM datasets for logistic regression.}
    \label{fig:loss}
\end{figure*}

\begin{figure*}[ht]
    \centering
    \subfigure[]{
    \includegraphics[width=0.20\textwidth]{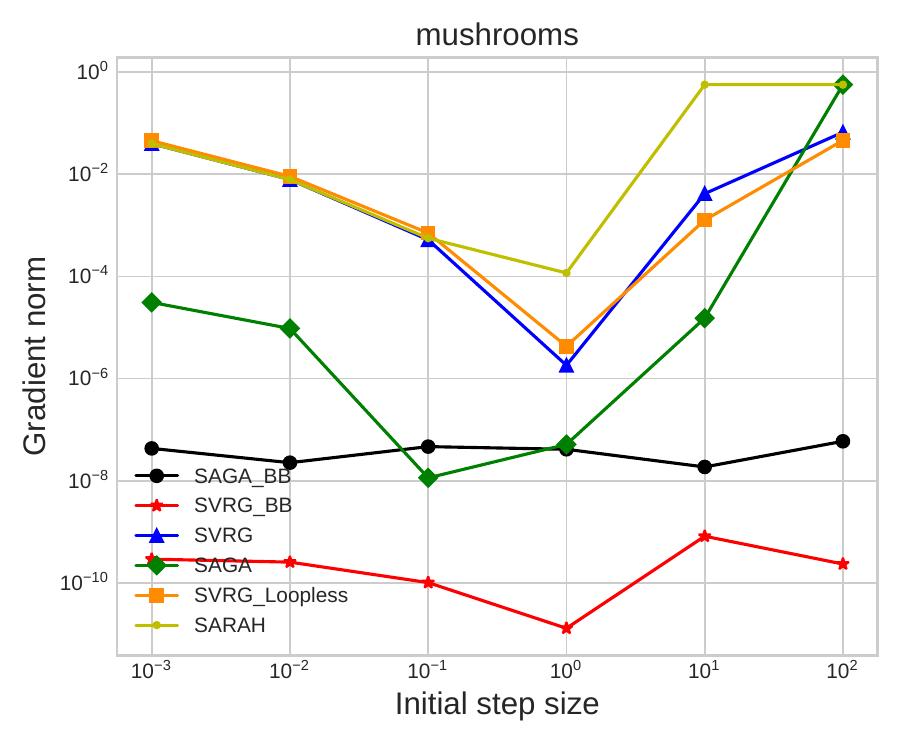}}
    \subfigure[]{
    \includegraphics[width=0.20\textwidth]{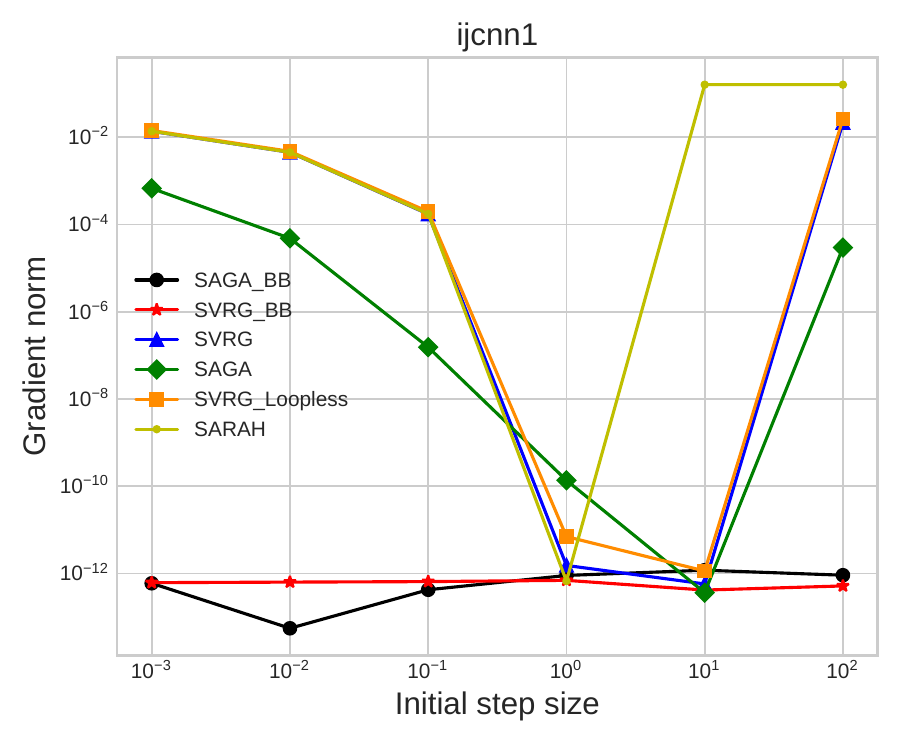}}
    \subfigure[]{
    \includegraphics[width=0.20\textwidth]{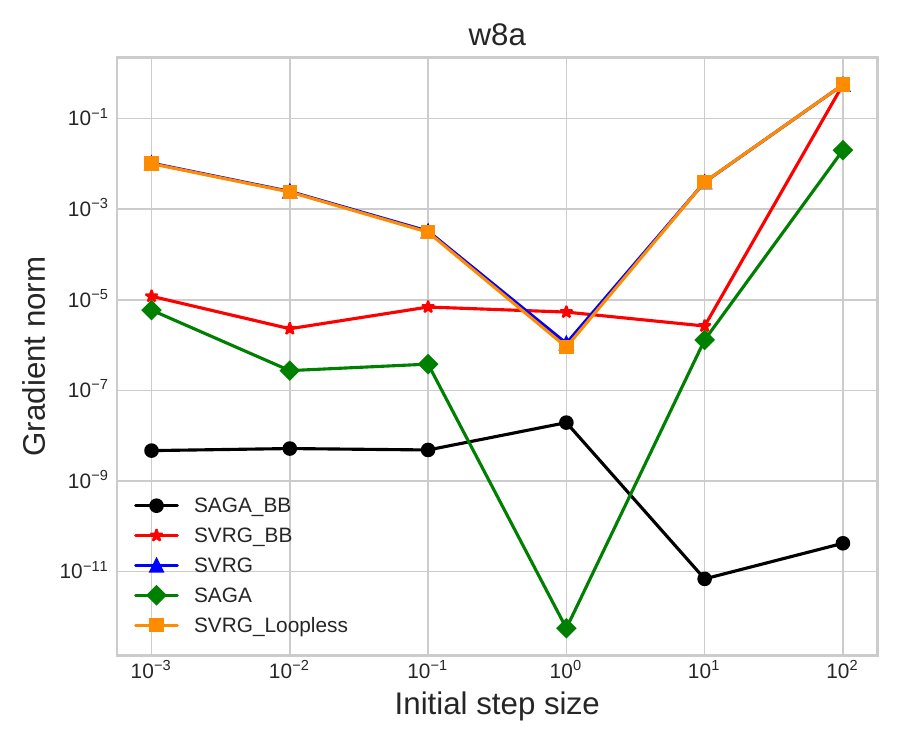}}
    \subfigure[]{
    \includegraphics[width=0.20\textwidth]{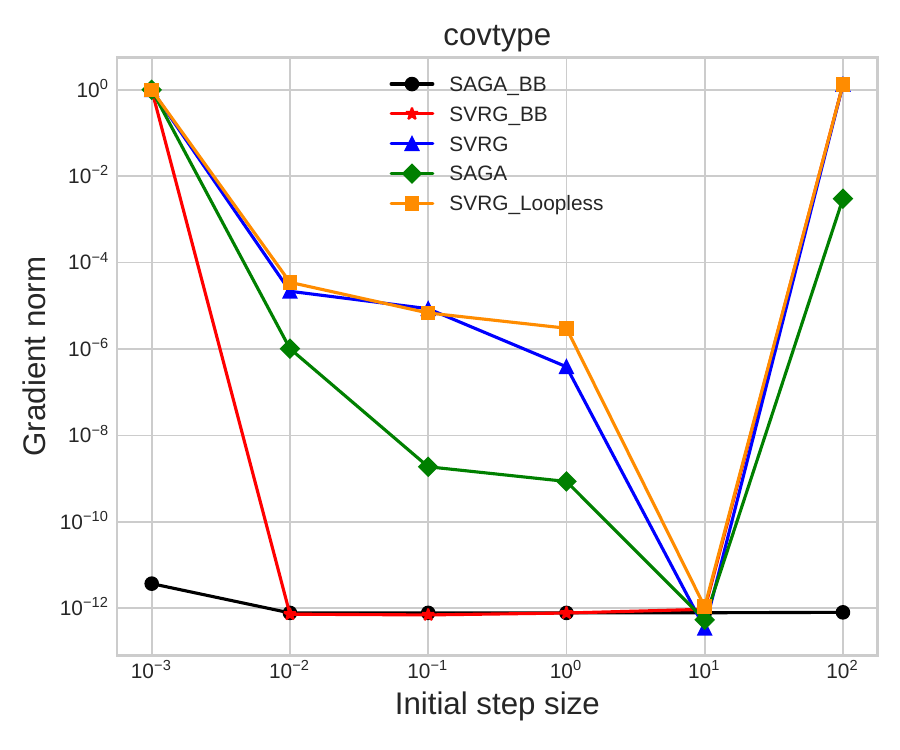}}
    \caption{Performance under different initial step-sizes and a fixed batch-size of 8. In some cases, we remove the curves so as not to clutter the plots.}
    \label{fig:grad_norm-step_size}
\end{figure*}

\section{Experiments}

In this section, we conduct extensive numerical experiments to verify the efficacy of SAGA-BB.
We apply our method and other variance-reduced methods for solving $\ell_2$-regularized logistic regression problems on binary classification tasks  with $a_i \in \mathbb{R}^d, b_i \in \{-1, +1 \} $, indexed by $i \in [n]$ : $f(x) = \frac{1}{n}\sum_{i=1}^{n}\log(1+\exp(-b_i a_i^T x))+\frac{\lambda}{2}\|x\|^2,$ where $(a_i, b_i)$ is the training sample. $\lambda$ is the regularization parameter which is set to $1/n$.

The datasets we choose are \emph{mushrooms, ijcnn1, w8a} and \emph{covtype}.\footnote{\mbox{All the datasets we used are available at https://www.csie} .ntu.edu.tw/˜cjlin/libsvmtools/datasets/.}
We compare SAGA-BB with SAGA, SVRG, SVRG-BB, SVRG-Loopless and SARAH.
For a fair comparison, we search over the initial step-size from $[10^{-3}, 10^{-2}, 10^{-1}, 1, 10, 100]$ and batch-size from $[1, 8, 16, 64]$ for each algorithm and each experiment.
As is common, we set $m=n/batch\_size$ for these methods.

We present our results in Figure \ref{fig:loss}, where we denote the average value by colored lines and standard deviation by filled areas, across 5 independent runs.
The vertical axis is the training objective minus the best loss, and the horizontal axis denotes the number of gradient evaluations normalized by $n$.
We compared all the methods with their best-tuning parameters.
We can see that SAGA-BB shows superior performance over other methods in all cases.
In addition, SAGA-BB significantly outperforms the original SAGA method, even on the dataset where SAGA performs well.

Figure \ref{fig:grad_norm-step_size} shows the final gradient norm with respect to the initial step-size,
presenting the robustness of SAGA-BB.
In some cases, the training objective minus the best loss can be zero at the final iteration.
Thus, we use the final gradient norm here.
We run all the methods with 120 effective gradient passes and different step-sizes.
To alleviate other parameters' impact, we choose a batch-size of 8, with which most of the other methods perform reasonably well.
In some cases, SARAH diverged, and we remove the curves accordingly so as not to clutter the plots.
We can see that most of the vanilla VR methods are sensitive to the initial step-size on some datasets.
For example, SARAH has a big oscillatory behavior on \emph{mushrooms}, and so does SVRG on \emph{covtype}.
Despite being equipped with BB step-size as well, SVRG-BB is inconsistent across the datasets.
% It is comparable with SAGA-BB on \emph{ijcnn} but shows an obvious oscillatory behaviour on \emph{w8a} dataset.
%
However, SAGA-BB shows a more consistent behavior that it always finds a good solution for all the datasets.
Thus, the experiment results demonstrate that SAGA-BB shows a strongly competitive performance.
Additional figures can be found in Appendix \ref{appendix: exps} for further confirming our claim.

\section{Conclusion}
In this paper, we propose SAGA-BB, a variant of SAGA with a novel adaptive step-size.
With a simple and clear construction, it exhibits a robust performance across different settings and random seeds.
SAGA-BB can automatically adjust the step-size by approximating the local geometry with stochastic information.
Theoretically, we prove that SAGA-BB keeps up with the oracle complexity achieved by SAGA.
For future work, it would be interesting to incorporate our adaptive step-size to other stochastic methods in various settings, such as the accelerated methods.

\bibliographystyle{named}
\bibliography{reference}

\begin{thebibliography}{}

\bibitem[\protect\citeauthoryear{Allen-Zhu and Orecchia}{2014}]{allen2014linear}
Zeyuan Allen-Zhu and Lorenzo Orecchia.
\newblock Linear coupling: An ultimate unification of gradient and mirror descent.
\newblock {\em arXiv preprint arXiv:1407.1537}, 2014.

\bibitem[\protect\citeauthoryear{Allen-Zhu}{2017}]{DBLP:conf/stoc/Zhu17}
Zeyuan Allen-Zhu.
\newblock Katyusha: the first direct acceleration of stochastic gradient methods.
\newblock In {\em Proceedings of the 49th Annual {ACM} {SIGACT} Symposium on Theory of Computing}, pages 1200--1205, 2017.

\bibitem[\protect\citeauthoryear{Armijo}{1966}]{pjm/1102995080}
Larry Armijo.
\newblock {Minimization of functions having Lipschitz continuous first partial derivatives.}
\newblock {\em Pacific Journal of Mathematics}, 16(1):1 -- 3, 1966.

\bibitem[\protect\citeauthoryear{Barzilai and Borwein}{1988}]{barzilai1988two}
Jonathan Barzilai and Jonathan~M Borwein.
\newblock Two-point step size gradient methods.
\newblock {\em IMA journal of numerical analysis}, 8(1):141--148, 1988.

\bibitem[\protect\citeauthoryear{Ben-Tal \bgroup \em et al.\egroup }{2001}]{BenTal2001TheOS}
A.~Ben-Tal, T.~Margalit, and A.~Nemirovski.
\newblock The ordered subsets mirror descent optimization method with applications to tomography.
\newblock {\em SIAM J. Optim.}, 12:79--108, 2001.

\bibitem[\protect\citeauthoryear{Burdakov \bgroup \em et al.\egroup }{2019}]{burdakov2019stabilized}
Oleg Burdakov, Yu-Hong Dai, and Na~Huang.
\newblock Stabilized barzilai-borwein method.
\newblock {\em arXiv preprint arXiv:1907.06409}, 2019.

\bibitem[\protect\citeauthoryear{Defazio \bgroup \em et al.\egroup }{2014}]{DBLP:conf/nips/DefazioBL14}
Aaron Defazio, Francis~R. Bach, and Simon Lacoste{-}Julien.
\newblock {SAGA:} {A} fast incremental gradient method with support for non-strongly convex composite objectives.
\newblock In {\em Advances in Neural Information Processing Systems}, pages 1646--1654, 2014.

\bibitem[\protect\citeauthoryear{Dubois{-}Taine \bgroup \em et al.\egroup }{2021}]{DBLP:journals/corr/abs-2102-09645}
Benjamin Dubois{-}Taine, Sharan Vaswani, Reza Babanezhad, Mark Schmidt, and Simon Lacoste{-}Julien.
\newblock {SVRG} meets adagrad: Painless variance reduction.
\newblock {\em CoRR}, abs/2102.09645, 2021.

\bibitem[\protect\citeauthoryear{Dvurechensky and Gasnikov}{2016}]{dvurechensky2016stochastic}
Pavel Dvurechensky and Alexander Gasnikov.
\newblock Stochastic intermediate gradient method for convex problems with stochastic inexact oracle.
\newblock {\em Journal of Optimization Theory and Applications}, 171(1):121--145, 2016.

\bibitem[\protect\citeauthoryear{Fletcher}{2005}]{fletcher2005barzilai}
Roger Fletcher.
\newblock On the barzilai-borwein method.
\newblock In {\em Optimization and control with applications}, pages 235--256. Springer, 2005.

\bibitem[\protect\citeauthoryear{Gorbunov \bgroup \em et al.\egroup }{2020}]{gorbunov2020stochastic}
Eduard Gorbunov, Marina Danilova, and Alexander Gasnikov.
\newblock Stochastic optimization with heavy-tailed noise via accelerated gradient clipping.
\newblock {\em Advances in Neural Information Processing Systems}, 33:15042--15053, 2020.

\bibitem[\protect\citeauthoryear{Johnson and Zhang}{2013}]{DBLP:conf/nips/Johnson013}
Rie Johnson and Tong Zhang.
\newblock Accelerating stochastic gradient descent using predictive variance reduction.
\newblock In {\em Advances in Neural Information Processing Systems}, pages 315--323, 2013.

\bibitem[\protect\citeauthoryear{Kovalev \bgroup \em et al.\egroup }{2020}]{DBLP:conf/alt/KovalevHR20}
Dmitry Kovalev, Samuel Horv{\'{a}}th, and Peter Richt{\'{a}}rik.
\newblock Don't jump through hoops and remove those loops: {SVRG} and katyusha are better without the outer loop.
\newblock In {\em Algorithmic Learning Theory}, 2020.

\bibitem[\protect\citeauthoryear{Li \bgroup \em et al.\egroup }{2020}]{DBLP:conf/icml/LiWG20}
Bingcong Li, Lingda Wang, and Georgios~B. Giannakis.
\newblock Almost tune-free variance reduction.
\newblock In {\em Proceedings of the 37th International Conference on Machine Learning}, 2020.

\bibitem[\protect\citeauthoryear{Liu \bgroup \em et al.\egroup }{2019}]{c:21}
Yan Liu, Congying Han, and Tiande Guo.
\newblock A class of stochastic variance reduced methods with an adaptive stepsize.
\newblock Forthcoming, 2019.

\bibitem[\protect\citeauthoryear{Loizou \bgroup \em et al.\egroup }{2021}]{DBLP:conf/aistats/LoizouVLL21}
Nicolas Loizou, Sharan Vaswani, Issam~Hadj Laradji, and Simon Lacoste{-}Julien.
\newblock Stochastic polyak step-size for {SGD:} an adaptive learning rate for fast convergence.
\newblock In {\em The 24th International Conference on Artificial Intelligence and Statistics}, 2021.

\bibitem[\protect\citeauthoryear{Malitsky and Mishchenko}{2020}]{DBLP:conf/icml/MalitskyM20}
Yura Malitsky and Konstantin Mishchenko.
\newblock Adaptive gradient descent without descent.
\newblock In {\em Proceedings of the 37th International Conference on Machine Learning}, 2020.

\bibitem[\protect\citeauthoryear{Nemirovski \bgroup \em et al.\egroup }{2009}]{DBLP:journals/siamjo/NemirovskiJLS09}
Arkadi Nemirovski, Anatoli~B. Juditsky, Guanghui Lan, and Alexander Shapiro.
\newblock Robust stochastic approximation approach to stochastic programming.
\newblock {\em {SIAM} J. Optim.}, 19(4):1574--1609, 2009.

\bibitem[\protect\citeauthoryear{Nemirovskij and Yudin}{1983}]{nemirovskij1983problem}
Arkadij~Semenovi{\v{c}} Nemirovskij and David~Borisovich Yudin.
\newblock Problem complexity and method efficiency in optimization.
\newblock 1983.

\bibitem[\protect\citeauthoryear{Nesterov}{2004}]{DBLP:books/sp/Nesterov04}
Yurii~E. Nesterov.
\newblock {\em Introductory Lectures on Convex Optimization - {A} Basic Course}, volume~87 of {\em Applied Optimization}.
\newblock Springer, 2004.

\bibitem[\protect\citeauthoryear{Nguyen \bgroup \em et al.\egroup }{2017}]{DBLP:conf/icml/NguyenLST17}
Lam~M. Nguyen, Jie Liu, Katya Scheinberg, and Martin Tak{\'{a}}c.
\newblock {SARAH:} {A} novel method for machine learning problems using stochastic recursive gradient.
\newblock In {\em Proceedings of the 34th International Conference on Machine Learning}, 2017.

\bibitem[\protect\citeauthoryear{Parikh and Boyd}{2014}]{parikh2014proximal}
Neal Parikh and Stephen Boyd.
\newblock Proximal algorithms.
\newblock {\em Foundations and Trends in optimization}, 1(3):127--239, 2014.

\bibitem[\protect\citeauthoryear{Pedregosa \bgroup \em et al.\egroup }{2011}]{scikit-learn}
F.~Pedregosa, G.~Varoquaux, A.~Gramfort, V.~Michel, B.~Thirion, O.~Grisel, M.~Blondel, P.~Prettenhofer, R.~Weiss, V.~Dubourg, J.~Vanderplas, A.~Passos, D.~Cournapeau, M.~Brucher, M.~Perrot, and E.~Duchesnay.
\newblock Scikit-learn: Machine learning in {P}ython.
\newblock {\em Journal of Machine Learning Research}, 12:2825--2830, 2011.

\bibitem[\protect\citeauthoryear{Robbins and Monro}{1951}]{10.1214/aoms/1177729586}
Herbert Robbins and Sutton Monro.
\newblock {A Stochastic Approximation Method}.
\newblock {\em The Annals of Mathematical Statistics}, 22(3):400 -- 407, 1951.

\bibitem[\protect\citeauthoryear{Roux \bgroup \em et al.\egroup }{2012}]{DBLP:conf/nips/RouxSB12}
Nicolas~Le Roux, Mark Schmidt, and Francis~R. Bach.
\newblock A stochastic gradient method with an exponential convergence rate for finite training sets.
\newblock In {\em Advances in Neural Information Processing Systems}, pages 2672--2680, 2012.

\bibitem[\protect\citeauthoryear{Shalev-Shwartz and Singer}{2007}]{shalev2007online}
Shai Shalev-Shwartz and Yoram Singer.
\newblock Online learning: Theory, algorithms, and applications.
\newblock 2007.

\bibitem[\protect\citeauthoryear{Shi \bgroup \em et al.\egroup }{2021}]{DBLP:journals/corr/abs-2102-09700}
Zheng Shi, Nicolas Loizou, Peter Richt{\'{a}}rik, and Martin Tak{\'{a}}c.
\newblock {AI-SARAH:} adaptive and implicit stochastic recursive gradient methods.
\newblock {\em CoRR}, abs/2102.09700, 2021.

\bibitem[\protect\citeauthoryear{Tan \bgroup \em et al.\egroup }{2016}]{DBLP:conf/nips/TanMDQ16}
Conghui Tan, Shiqian Ma, Yu{-}Hong Dai, and Yuqiu Qian.
\newblock Barzilai-borwein step size for stochastic gradient descent.
\newblock In {\em Advances in Neural Information Processing Systems}, pages 685--693, 2016.

\bibitem[\protect\citeauthoryear{Vrahatis \bgroup \em et al.\egroup }{2000}]{vrahatis2000class}
Michael~N Vrahatis, George~S Androulakis, JN~Lambrinos, and George~D Magoulas.
\newblock A class of gradient unconstrained minimization algorithms with adaptive stepsize.
\newblock {\em Journal of Computational and Applied Mathematics}, 114(2):367--386, 2000.

\bibitem[\protect\citeauthoryear{Zhou \bgroup \em et al.\egroup }{2019}]{zhou2019direct}
Kaiwen Zhou, Qinghua Ding, Fanhua Shang, James Cheng, Danli Li, and Zhi-Quan Luo.
\newblock Direct acceleration of saga using sampled negative momentum.
\newblock In {\em The 22nd International Conference on Artificial Intelligence and Statistics}, pages 1602--1610, 2019.

\bibitem[\protect\citeauthoryear{Zhou \bgroup \em et al.\egroup }{2020}]{zhou2020amortized}
Kaiwen Zhou, Yanghua Jin, Qinghua Ding, and James Cheng.
\newblock Amortized nesterov’s momentum: A robust momentum and its application to deep learning.
\newblock In {\em Conference on Uncertainty in Artificial Intelligence}, pages 211--220, 2020.

\end{thebibliography}
\appendix
\onecolumn
\section{Implementation of SAGA-BB}
\label{appendix: implement}
The memory-efficient implementation of SAGA-BB at iteration $k$ is presented in Algorithm \ref{alg:impt}.
It is inspired by the fact that the objectives $f_i(x)$ in many ML tasks can be interpreted as a function of scalars \citep{zhou2019direct}, i.e., $\psi_i(\left <a_i,x \right>).$
We only store the scalar $\left <a_i, \phi_i^k \right>$ instead of the vector $\nabla f_i(\phi_i^k)$ for each $i$.
Based on this, we can have the value of $s_k^Ty_k$ and $y^T_ky_k$ without extra storage requirement as shown in Algorithm \ref{alg:impt}, whereas other adaptive step-sizes need to store the points vector $\phi_i^k$.
Thus, we achieve the $\mathcal{O}(n)$ storage requirement for SAGA-BB in practice.

\begin{algorithm}[tb]
\caption{Implementation of SAGA-BB}
\label{alg:impt}
\textbf{Stored}: Scalar table $\Phi^k$ where $\Phi^k_i = \left <a_i,\phi_i^k \right>$, and $\mu^k.$
\begin{algorithmic}[1] %[1] enables line numbers
\STATE \textbf{At iteration $k$}
\IF {$k \mod m == 0$}
\STATE
$ \eta_k = \frac{\sum_{i=k-m}^{k} s_{i}^Ty_{i}}{\alpha \sum_{i=k-m}^{k} y_{i}^Ty_i}$
\ENDIF
\STATE Randomly choose a sample $i$ and compute the gradient estimator with the scalar table. Store $\Phi^{k+1}_i=\left<a_i,x_k\right>$ with other entries unchanged.
\STATE 
$\widetilde{\nabla}_k = (\nabla \psi_i(\Phi^{k+1}_i) -\nabla \psi_i(\Phi^k_i)) \cdot d_i + \mu^k$
\STATE
$\mu_{k+1}=\mu_{k} + \frac{1}{n}(\nabla \psi_i(\Phi^{k+1}_i) -\nabla \psi_i(\Phi^k_i)) \cdot d_i$
\STATE
$s^T_k y_k = \Phi^{k+1}_i \cdot \nabla \psi_i(\Phi^{k+1}_i)  + \Phi^{k+1}_i \cdot \nabla \psi_i(\Phi^{k}_i) + \Phi^k_i \cdot \nabla \psi_i(\Phi^{k+1}_i) + \Phi^k_i \cdot \nabla \psi_i(\Phi^k_i)$
\STATE
$y_k^T y_k = \|\nabla \psi_i(\Phi^{k+1}_i) \cdot d_i \|^2 + 2(\nabla \Phi^{k+1}_i \cdot \nabla \psi_i(\Phi^{k}_i))\|d_i\|^2 + \|\nabla \psi_i(\Phi^{k}_i) \cdot d_i \|^2$
\STATE
$x_{k+1}=\arg \min_x \left\{\left\langle \widetilde{\nabla}_{k}, x \right\rangle  + \frac{1}{2\eta_k}\|x-x_k\|^2+ h(x)\right\}$
\end{algorithmic}
\end{algorithm}

\section{ALGORITHM FOR MINI-BATCH}
\label{appendix: mini-batch}

Let $|B|$ denote the batch size.
At each iteration $k$, we will pick a mini-batch $B$ from ${1,..,n}$ of size $|B|$ and the mini-batch gradients can be represented as:
\begin{equation}
    \nabla f_B(\phi_B^k) = \frac{1}{|B|}\sum_{j \in B}f_j(\phi_j^k).
\end{equation}
We also use $\phi_B$ to denote the the set of points $\phi_j$, for each $j \in B$.

\begin{algorithm}[h!]
\caption{SAGA-BB}
\label{alg:algorithm}
\textbf{Input}: Max epochs $K$, batch-size $|B|$, initial point $x_1$, initial step size $\eta_1$, update frequency $m$.  \\
\textbf{Initialize}: Initial the gradient table $\nabla f_i(\phi_i^{1})$ table with $\phi_i^{1}=x_{1}$, and the average $\mu = \frac{1}{n}\sum_{i=1}^n \nabla f_i(\phi_i^{1})$.  \\

\begin{algorithmic}[1] %[1] enables line numbers
\FOR{$k = 1,2,\dots, K$}
\IF {$k \mod m == 0$}
\STATE
$ \eta_k = \frac{\sum_{i=k-m}^{k} s_{i}^Ty_{i}}{\alpha \sum_{i=k-m}^{k} y_{i}^Ty_i}$
\STATE
\ENDIF

\STATE Randomly pick a batch $B$, take $\phi_j^{k+1}=x_{k}$, store $\nabla f_{j}(\phi_{j}^{k+1})$ with other entries unchanged for every $j \in B$.
\STATE
$\widetilde{\nabla}_{k} =\nabla f_{B}(\phi_{B}^{k+1})-\nabla f_{B}(\phi_{B}^{k})+\mu_{k}$
\STATE
$y_{k}=\sum_{j \in B}(\nabla f_{j}(\phi_{j}^{k+1})-\nabla f_{j}(\phi_{j}^{k}))$
\STATE
$s_{k}=\sum_{j \in B}(\phi_{j}^{k+1}-\phi_{j}^{k})$
\STATE
$\mu_{k+1}=\mu_{k} + \frac{1}{n}\sum_{j \in B}(\phi_{j}^{k+1}-\phi_{j}^{k})$
\STATE
%xk+1=xk−ηk˜∇kx_{k+1} = x_{k} - \eta_k \widetilde{\nabla}_{k}
$x_{k+1}=\arg min_x \left\{\left\langle \eta_k \widetilde{\nabla}_{k}, x \right\rangle  + V_d(x,x_k)+\eta_k h(x)\right\}$
\STATE

\ENDFOR
\STATE \textbf{return} $x_{K+1}$
\end{algorithmic}
\end{algorithm}

\begin{algorithm}[!]
	\caption{Restarted SAGA-BB (R-SAGA-BB)}
	\label{Restart}
	\textbf{Input}: Max epochs $\tau$, number of iterations $K$ of SAGA-BB, initial point $x_1$, initial step size $\eta_1$.  \\
    \textbf{Initialize}: Initial the gradient table $\nabla f_i(\phi_i^{1})$ table with $\phi_i^{1}=x_{1}$, and the average $\mu = \frac{1}{n}\sum_{i=1}^n \nabla f_i(\phi_i^{1})$.  \\

	\begin{algorithmic}[1] %[1] enables line numbers
		%\STATE Let t=0t=0.
		\FOR{$t = 1,2,\dots, \tau$}
		\STATE 		 Run SAGA-BB (Algorithm \ref{alg:algorithm}) for $K$ iterations with batch-size $|B|$, step size $\gamma$, starting point $\bar{x}_t,$ and update frequency $m$. Define the output of SAGA-BB by $\bar{x}_{t+1}.$ 
		\ENDFOR
		\STATE \textbf{return} $\bar{x}_{\tau+1}.$
	\end{algorithmic}
\end{algorithm}

\section{Importance of Non-Euclidean Norms}
\label{appendix: importance}
Euclidean norm space is common and widely considered in both academia and industry.
However, there exist some problems related to non-Euclidean norms that are also crucial to study both in theory \citep{allen2014linear} and industry (e.g., Positron Emission Tomography (PTE) can be seen as a Maximum Likelihood Estimate problem related to $\ell_1$ norm case over the large simplex \citep{BenTal2001TheOS}).
Katyusha proposed by \cite{DBLP:conf/stoc/Zhu17}, as a variant of SVRG, can be extended to non-Euclidean norm smoothness settings.
However, there is still a void in the existing analysis for SAGA and its variants.
The main reason is that the theory of SAGA crucially relies on the $\ell_2$ norm space.
%is hardly used in the general norm space since it is well designed in  ℓ2\ell_2 norm space
%\WM{If the theory is well designed/developed, why is there still a void?}{\color{blue} Thanks, the revised one is above. Just showing the theory of original SAGA only work in euclidean norm space and details are as follows.} In particular, \citet{DBLP:conf/nips/DefazioBL14} defined a Lyapunov function as an ``intermediary '' which is critical to construct a contraction with some key lemmas only working in ℓ2\ell_2 norm space.
%
To tackle this problem, we propose a new Lyapunov function which supports arbitrary norms and show that it can also derive a characteristic inequality.
In the following subsections, we conduct a comprehensive study on the convergence analysis.

\section{PROOF OF LEMMA \ref{VB}}
\begin{align*}
	\mathbb{E}\left[\left\|\nabla f(x_k)-\widetilde{\nabla}_{k}\right\|^2_*\right]&\overset{(a)}{\leq}  \mathbb{E}\left[\left\|\nabla f_B(x_k)-\nabla f_B(\phi_B)\right\|^2_*\right]\\
	&{}\overset{(b)}{\leq} 2\mathbb{E}\left[\left\|\nabla f_B(x_k)-\nabla f_B(x^*)\right\|^2_*\right]+2\mathbb{E}\left[\left\|\nabla f_B(\phi_B)-\nabla f_B(x^*)\right\|^2_*\right]\\
	&{}\overset{(c)}{\leq} 4L\left[f(x_k)-f(x^*)- \langle \nabla f(x^*), x_k - x^*\rangle \right]\\
	&{} + 4L\left[\frac{1}{n}\sum_{i=1}^{n}f_i(\phi_i^k)-f(x^*)- \frac{1}{n}\sum_{i=1}^{n}\langle \nabla f_i(x^*), \phi_i^k - x^*\rangle \right], 
	\end{align*}
	where $(a)$ uses $\mathbb{E}[\|X-\mathbb{E}X\|^2_*] \leq \mathbb{E}[\|X\|^2_*],$ $(b)$ uses $\|a+b\|^2_* \leq 2\|a\|^2_* + 2\|b\|^2_*$ and $(c)$ follows from Theorem 2.1.5 in \cite{DBLP:books/sp/Nesterov04}.

\section{PROOF OF THEOREM \ref{theo1}}
\label{appendix: theorem2}
In this case, we modify the Lyapunov function as 
$$T_k:=\frac{1}{\omega}\left[\frac{1}{n}\sum_{i=1}^{n}F_i(\phi_i^k)- F(x^*)-\frac{1}{n}\sum_{i=1}^{n}\langle \nabla F_i(x^*),\phi_i^k - x^* \rangle \right].$$
Similarly, we start with per-iteration analysis in mini-batch setting. Using convexity of $f(\cdot)$ at $x_k,$ we have
	\begin{equation} \label{first}
		f(x_k)-f(x^*) \leq \langle\nabla f(x_k)-\widetilde{\nabla}_{k}, x_k -x^*\rangle+\langle \widetilde{\nabla}_{k}, x_k-x_{k+1}\rangle+\langle \widetilde{\nabla}_{k}, x_{k+1}-x^*\rangle.
	\end{equation}
	Taking expectation with respect to sample batch size $|B|$, we obtain 
	\begin{equation} \label{1}
		f(x_k)- f(x^*) \leq  \mathbb{E}\left[\langle\widetilde{\nabla}_{k}, x_k-x_{k+1}\rangle\right] +\mathbb{E}\left[\langle\widetilde{\nabla}_{k}, x_{k+1}- x^*\rangle\right].	
	\end{equation}
	It is direct to upper bound the first term on the right side using smoothness argument, i.e.,
	$$f(x_{k+1})- f(x_k) \leq \langle\nabla f(x_k)-\widetilde{\nabla}_{k}, x_{k+1} -x_k\rangle+\langle \widetilde{\nabla}_{k}, x_{k+1}-x_k\rangle+\frac{L}{2}\|x_{k+1}-x_k\|^2.$$
	Taking expectation with respect to sample $k$ and re-arranging, we obtain
	\begin{equation}\label{2}
	    \begin{split}
	    \mathbb{E}\left[\langle\widetilde{\nabla}_{k}, x_k-x_{k+1}\rangle\right] \leq{}&f(x_k)-\mathbb{E}[f(x_{k+1})] + \mathbb{E}\left[\langle\nabla f(x_k)-\widetilde{\nabla}_{k}, x_{k+1} -x_k\rangle\right]\\
	    &{}+\frac{L}{2}\mathbb{E}\left[\|x_{k+1}-x_k\|^2\right].
	    \end{split}
	\end{equation}
	Note that
	\begin{equation}\label{4}
		\begin{split}
		&\mathbb{E}[T_{k+1}]\\
		={}&\frac{1}{\omega}\bigg[(1-\frac{|B|}{n})\frac{1}{n}\sum_{i=1}^{n}F_i(\phi_i^k)- F(x^*) 
        - (1-\frac{|B|}{n})\frac{1}{n}\sum_{i=1}^{n}\langle \nabla F_i(x^*),\phi_i^k - x^* \rangle-\frac{|B|}{n}\langle\nabla F(x^*), x_k-x^*\rangle\\
	    {}&+\frac{|B|}{n}F(x_k) \bigg]\\
		\overset{(d)}{=}{}& T_k-\frac{|B|}{n}T_k+\frac{|B|}{n\omega}(F(x_k)-F(x^*)).%+\frac{|B|c}{n}V_d(x^*,x_k)-c(V_d(x^*,x_k)-\mathbb{E}[V_d(x^*,x_{k+1})]).
		\end{split}
	\end{equation}
	where $(d)$ using the optimality of $x^*$ that $\nabla F(x^*) = 0.$
	Upper-bounding (\ref{1}) using (\ref{2}), Lemma \ref{SCB}, (\ref{4}) and re-arranging, we obtain
	\begin{equation} \label{core}
	\begin{split}
		&\mathbb{E}[T_{k+1}]+\mathbb{E}[f(x_{k+1})]-f(x^*)\\ \leq{}& \mathbb{E}\left[\langle\nabla f(x_k)-\widetilde{\nabla}_{k}, x_{k+1} -x_k\rangle\right]+\frac{L}{2}\mathbb{E}\left[\|x_{k+1}-x_k\|^2\right]+\mathbb{E}[T_{k+1}]\\
		{}&+ \frac{1}{\eta_k}(V_d(x^*,x_{k+1})-\mathbb{E}[V_d(x^*,x_{k+1})]-\mathbb{E}[V_d(x_{k+1},x_k)])+h(x^*)-\mathbb{E}[h(x_{k+1})]\\
		\overset{(a)}{\leq}{}& \frac{1}{2\beta}\mathbb{E}\left[\left\|\nabla f(x_k)-\widetilde{\nabla}_{k}\right\|^2_*\right]+\left(\frac{\beta}{2}+\frac{L}{2}-\frac{1}{2\eta_k}\right)\mathbb{E}\left[\|x_{k+1}-x_k\|^2\right]\\
		{}&+\frac{1}{\eta_k}(V_d(x^*,x_{k})-\mathbb{E}[V_d(x^*,x_{k+1})])+\mathbb{E}[T_{k+1}]+h(x^*)-\mathbb{E}[h(x_{k+1})]\\
		\overset{(b)}{\leq}{}&T_k+\left(\frac{2L\omega}{\beta}-\frac{|B|}{n}\right)T_k + \left(\frac{2L}{\beta}+\frac{|B|}{n\omega}\right)[F(x_k)-F(x^*)]\\
		{}&+ \left(\frac{\beta}{2}+\frac{L}{2}-\frac{1}{2\eta_k}\right)\mathbb{E}\left[\|x_{k+1}-x_k\|^2\right]+ \frac{1}{\eta_k}(V_d(x^*,x_{k})-\mathbb{E}[V_d(x^*,x_{k+1})])\\
		{}&- \frac{2L}{\beta}\left[\frac{1}{n}\sum_{i=1}^{n}h(\phi_i^k)- h(x^*)-\frac{1}{n}\sum_{i=1}^{n}\langle \partial h(x^*),\phi_i^k - x^* \rangle \right]+h(x^*)-\mathbb{E}[h(x_{k+1})]\\
		{}&-\frac{2L}{\beta}[h(x_k)-h(x^*)-\langle \partial h(x^*),x_k-x^* \rangle],
		\end{split}
	\end{equation}
	where (a) uses Young's inequality with $\beta > 0 $ and that $d(\cdot)$ is 1-strongly convex, (b) follows from Lemma \ref{VB} and (\ref{4}).
	Letting $\beta = 8L,\omega = \frac{4|B|}{n}$ and $9L \leq \frac{1}{\eta_k} \leq \alpha L$, we obtain 
	\begin{align*}
		&\mathbb{E}[T_{k+1}]+\mathbb{E}[F(x_{k+1})]-F(x^*)\\ \leq{}& T_k + \frac{1}{2}[F(x_k)-F(x^*)]+ \frac{1}{\eta_k}(V_d(x^*,x_{k})-\mathbb{E}[V_d(x^*,x_{k+1})]).
	\end{align*}
% 	Note that ηkdef=∑ki=k−msTiyiα∑ki=k−myTiyi, \eta_k \overset{\text{def}}{=}   \frac{\sum_{i=k-m}^{k} s_{i}^Ty_{i}}{\alpha \sum_{i=k-m}^{k} y_{i}^Ty_i}, we can set the constant α≥9\alpha \geq 9, i.e., \[\eta_k=\begin{cases}
% \frac{\sum_{i=k-m}^{k} s_{i}^Ty_{i}}{\alpha \sum_{i=k-m}^{k} y_{i}^Ty_i}, &\text{if} \frac{\sum_{i=k-m}^{k} s_{i}^Ty_{i}}{\alpha \sum_{i=k-m}^{k} y_{i}^Ty_i} \leq \frac{1}{9L},\\
% \frac{1}{9L},&\text{otherwise}.\end{cases}\]
	Thus, 
    \begin{equation} \label{convex contraction}
	\begin{split}
	&\frac{1}{2}[\mathbb{E}[F(x_{k+1})]-F(x^*)]\\ \leq{}& T_k-\mathbb{E}[T_{k+1}] + \frac{1}{2}[F(x_k)-\mathbb{E}\left[F(x_{k+1})\right]]+ \frac{1}{\eta_k}(V_d(x^*,x_{k})-\mathbb{E}[V_d(x^*,x_{k+1})])
	\end{split}
    \end{equation}
	Summing the above inequality from $k=1,\cdots, K$ and taking expectation with respect to all randomness, we obtain
	$$\frac{1}{2}\sum_{i=1}^{K}\left[\mathbb{E}\left[F(x_{k+1})\right]\right]-KF(x^*) \leq T_1+\frac{1}{2}F(x_1)+\frac{1}{\eta_k}V_d(x^*,x_1)$$
	Using Jensen's inequality on the left side, we obtain the output after $K$ iterations of SAGA-BB:
	$$\mathbb{E}[F(\bar{x}_{K+1})] - F(x^*) \leq \frac{2}{K}\left[\frac{n}{4|B|}[F(x_1)-F(x^*)]+\frac{1}{2}F(x_1)+\alpha LV_d (x^*,x_1)\right].$$
Therefore, with $|B|=1$ we complete the proof of Theorem \ref{theo1}.

\section{PROOF OF THEOREM \ref{theo2}}
    Consider the first run of SAGA-BB (Algorithm \ref{alg:algorithm}), i.e., $t =0.$ Recall the proof of Theorem 1, we sum the inequality (\ref{convex contraction}) from $k=1,\cdots, K$ and taking expectation with respect to all randomness, we obtain
    \begin{align*}
        	\frac{1}{2}\sum_{i=1}^{K}\left[\mathbb{E}\left[F(x_{k+1})\right]\right]-KF(x^*) \leq{} & T_1+\frac{1}{2}[F(x_1)-\mathbb{E}[F(x_{k+1})]]+\frac{1}{\eta_k}V_d(x^*,x_1) \\
         \leq{}& T_1+\frac{1}{2}[F(x_1)-F(x^*)]+\frac{1}{\eta_k}V_d(x^*,x_1).
    \end{align*} 
    Then by using Jensen's inequality on the left side, we obtain the output after $K$ iterations of SAGA-BB:
	$$\mathbb{E}[F(\bar{x}_{K})] - F(x^*) \leq \frac{2}{K}\left[\left(\frac{n}{4|B|}+\frac{1}{2}\right)(F(x_1)-F(x^*))+\alpha LV_d (x^*,x_1)\right].$$
    Observed
	that the proof of Theorem 1 still holds if we substitute $V_d(x^*, x_1)$ everywhere by its upper bound $\frac{1}{\mu}(F(x_1) - F(x^*)).$ 
    By using generalized strong convexity of $f(\cdot),$ we have
    \begin{align*}
        \mu V_d(x^*, x_1) \leq{}& f(x_1)-f(x^*)-  \langle \nabla f(x^*),x_1-x^* \rangle \\
        \leq{}& f(x_1)-f(x^*)-  \langle \nabla f(x^*),x_1-x^* \rangle + h(x_1) - h(x^*) -\langle  \partial h(x^*),  x_1 - x^*\rangle\\
        \leq{}& F(x_1) - F(x^*) - \langle \nabla F(x^*), x_1 - x^* \rangle\\
        \leq{}& F(x_1) - F(x^*)
    \end{align*}
where the second inequality uses the convexity of $h(\cdot)$ and $\partial h(\cdot)$ denotes its sub-gradient. It implies that when $t=1,$ i.e., after $K$ iterations of SAGA-BB, we have
    $$\mathbb{E}[F(\bar{x}_{1})] - F(x^*) \leq \frac{2}{K}\left[\left(\frac{n}{4|B|}+\frac{1}{2}+\frac{\alpha L }{\mu}\right)(F(x_1)-F(x^*))\right].$$
Thus taking $K=\frac{\mu n + 6\alpha |B| L}{\mu |B|},$ we have 
	$$\mathbb{E}[F(\bar{x}_{1})] - F(x^*) \leq \frac{1}{2}(F(x_1)-F(x^*)).$$
Then, by induction, we can show that for arbitrary $t=1,2,\ldots,\tau$, the inequality 
    $$\mathbb{E}[F(\bar{x}_{\tau+1})] - F(x^*) \leq \frac{1}{2^{\tau}}(F(x_1)-F(x^*)).$$
Therefore, with $|B|=1$ we complete the proof of Theorem \ref{theo2}.

\section{ADDITIONAL EXPERIMENTS}
\label{appendix: exps}
Following the experiments presented in the main paper, we further evaluate the performance of our methods with additional experiments.

\subsection{Additional experiments with stable methods}\label{sub: stable methods}
\begin{figure}[ht]
    \centering
    \subfigure[]{
    \includegraphics[width=0.32\textwidth]{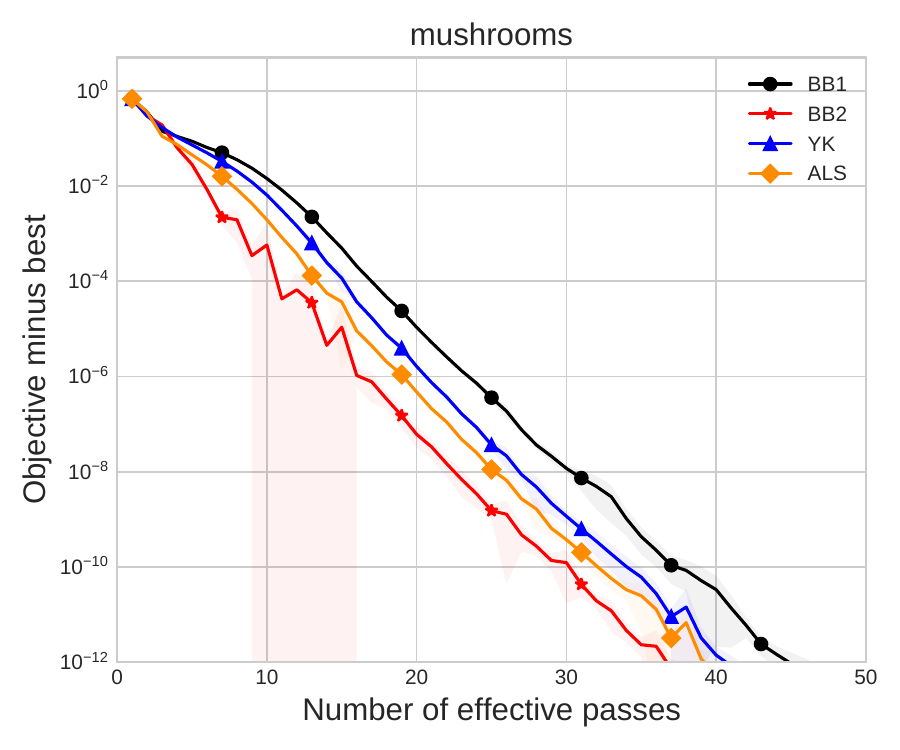}}
    \subfigure[]{
    \includegraphics[width=0.32\textwidth]{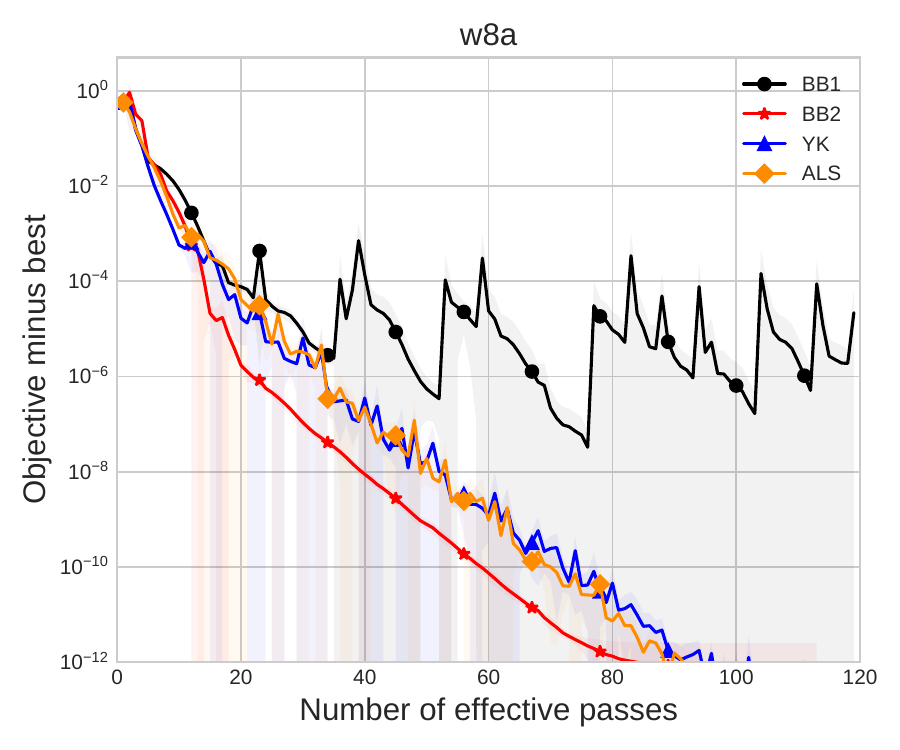}}
    \subfigure[]{
    \includegraphics[width=0.32\textwidth]{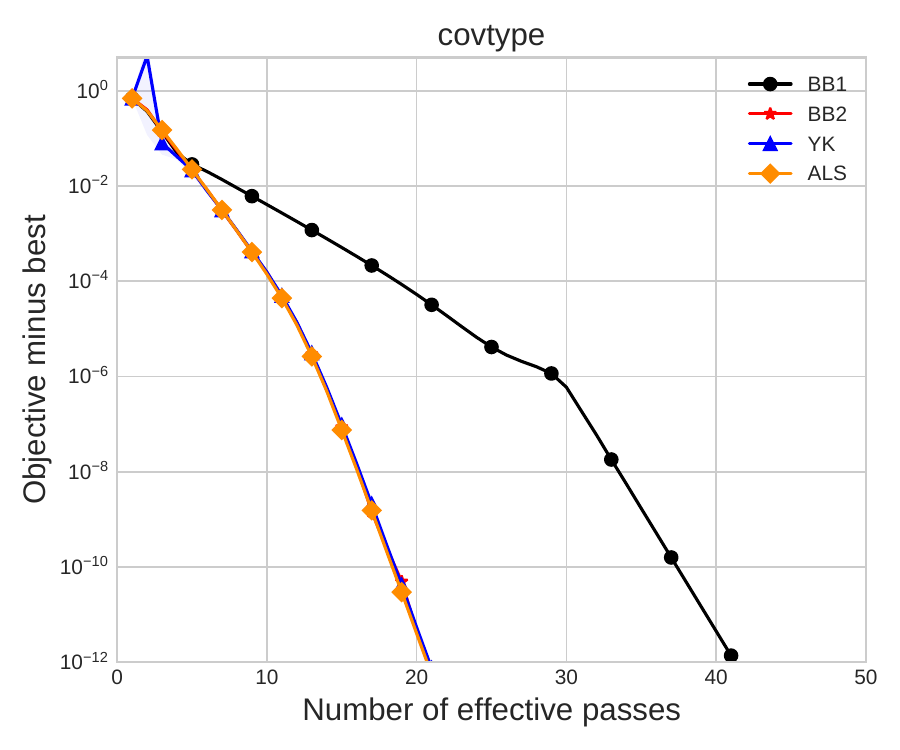}}
    \caption{Performance of different stable methods on different datasets.}
    \label{fig: additional stable methods}
\end{figure}

In this subsection, we evaluate the stable methods on different datasets, presented in Figure \ref{fig: additional stable methods}. Among these stable methods, BB2 shows a more consistent performance across different datasets.

\subsection{Additional experiments for logistic loss}
This subsection presents additional numerical results to present the robustness and the competitive performance of SAGA-BB.
\begin{figure}[ht]
    \centering
    \subfigure[]{
    \includegraphics[width=0.23\textwidth]{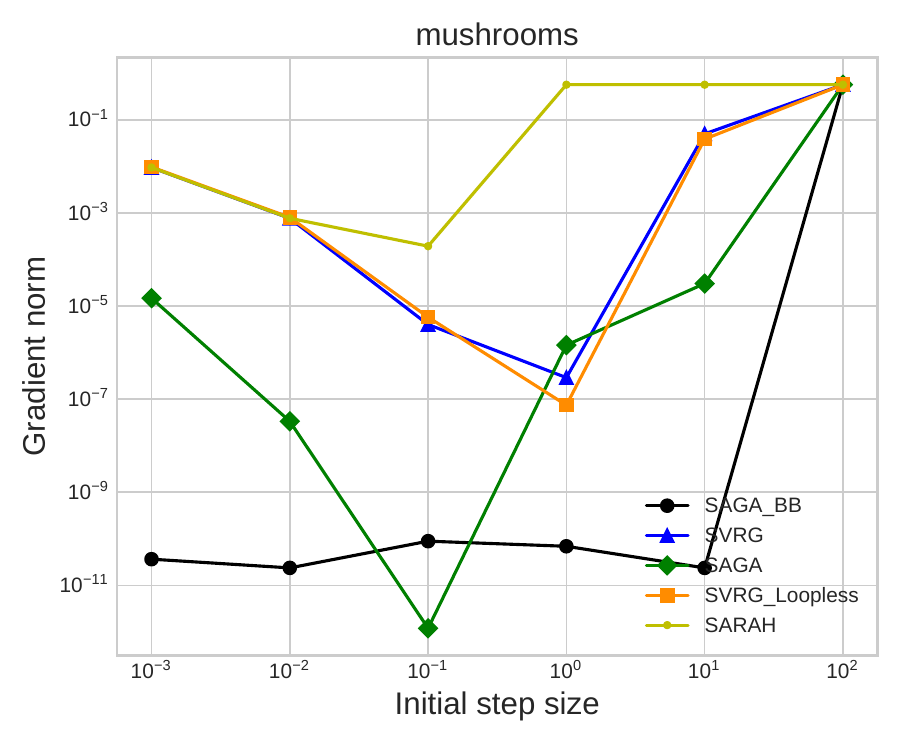}}
    \subfigure[]{
    \includegraphics[width=0.23\textwidth]{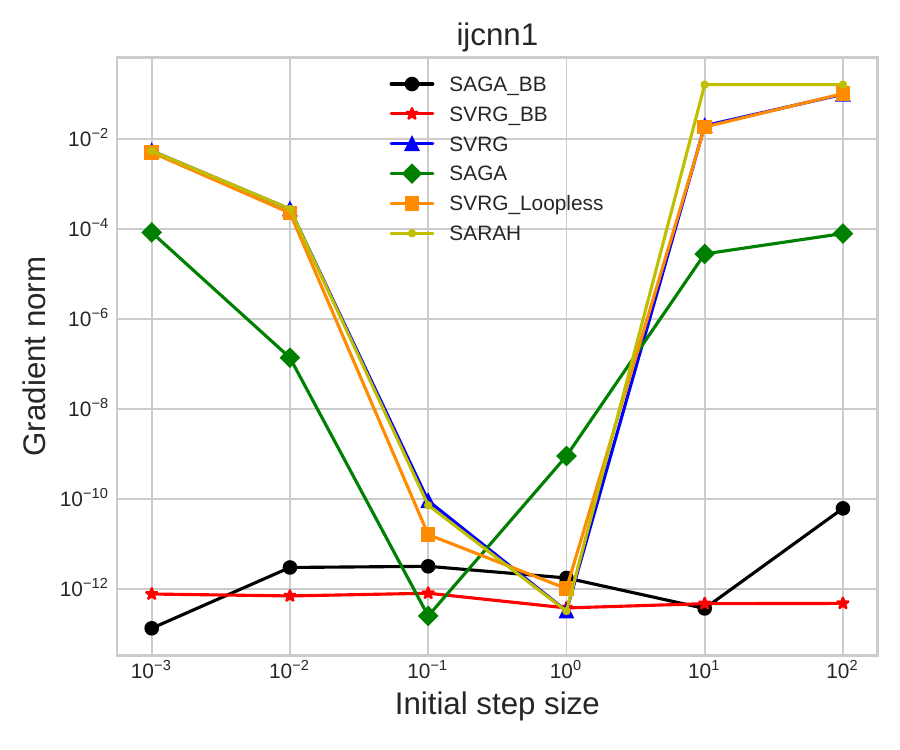}}
    \subfigure[]{
    \includegraphics[width=0.23\textwidth]{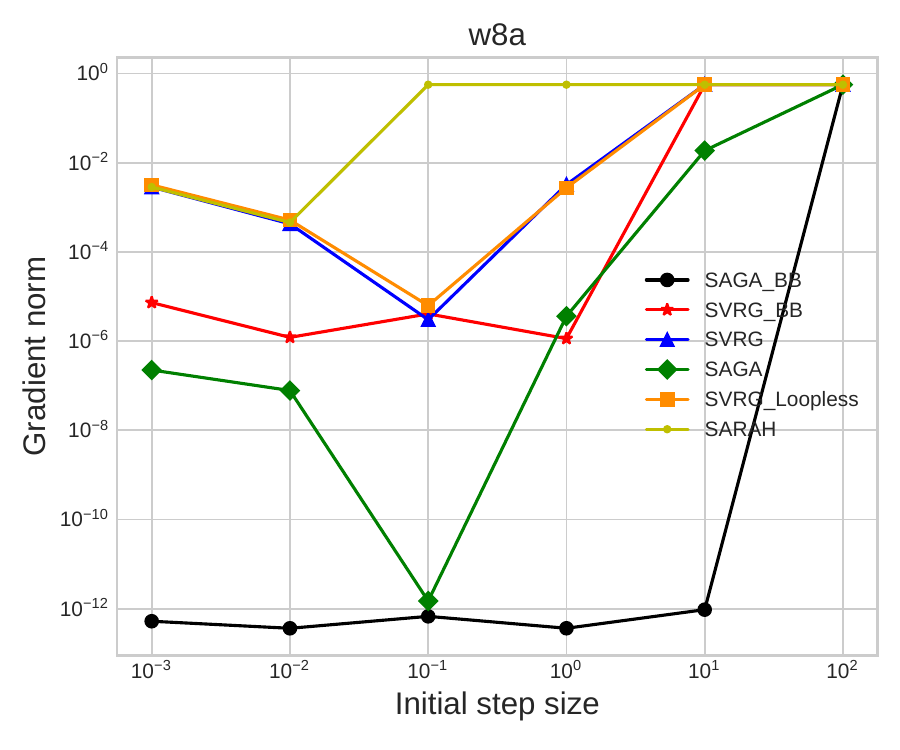}}
    \subfigure[]{
    \includegraphics[width=0.23\textwidth]{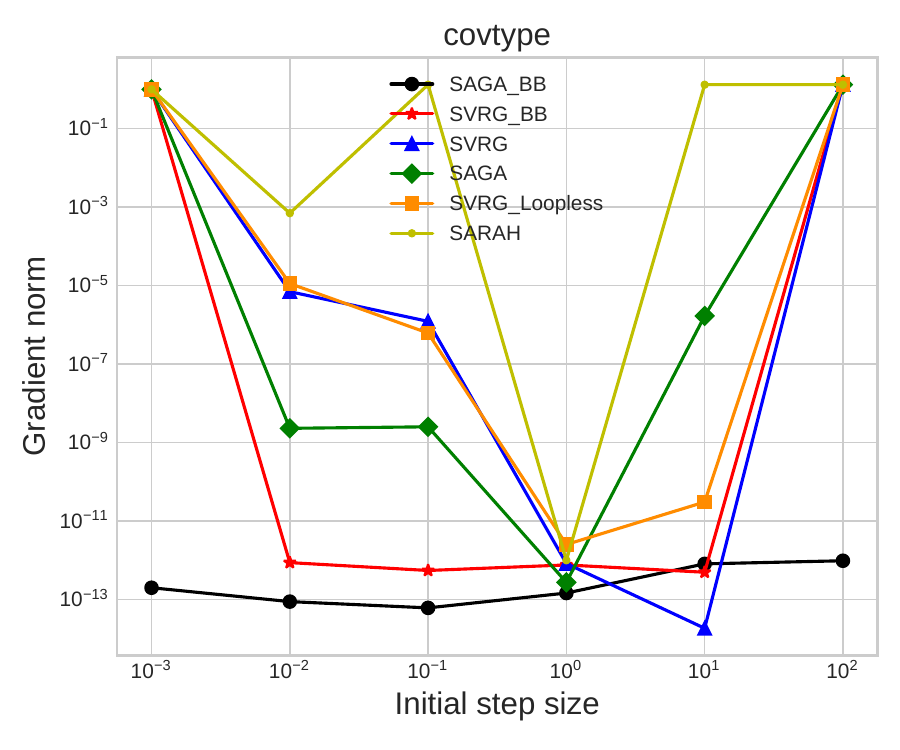}}
    \caption{Performance under different initial step size and fixed batch-size of 1. For the plot to be readable, we limit the gradient norm to a maximum value of 1.}
    \label{fig: logistic loss bz1}
\end{figure}

\begin{figure}[ht]
    \centering
    \subfigure[]{
    \includegraphics[width=0.23\textwidth]{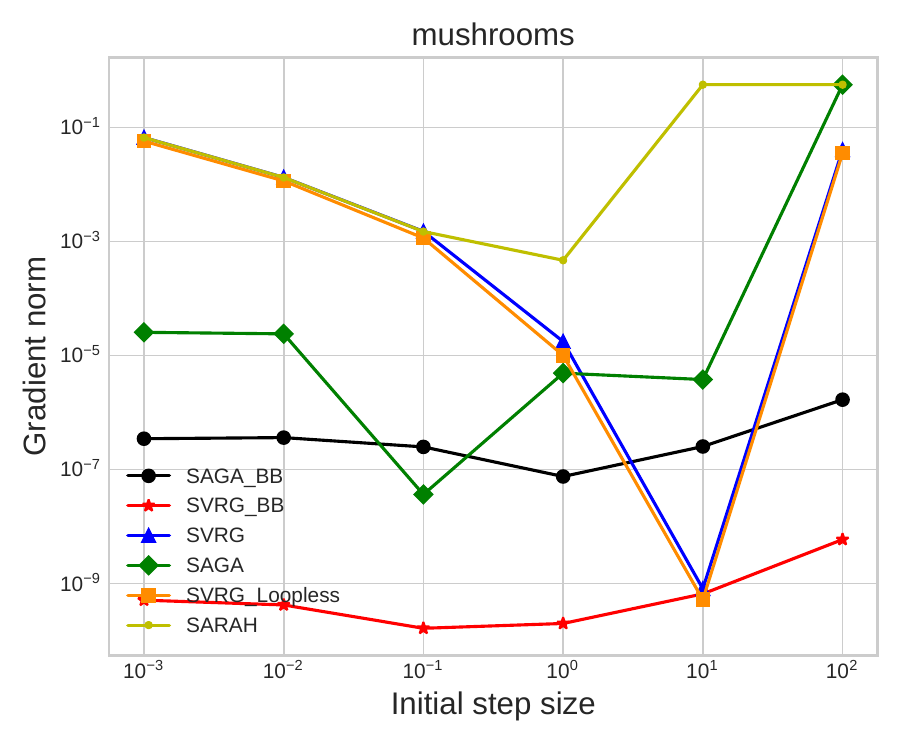}}
    \subfigure[]{
    \includegraphics[width=0.23\textwidth]{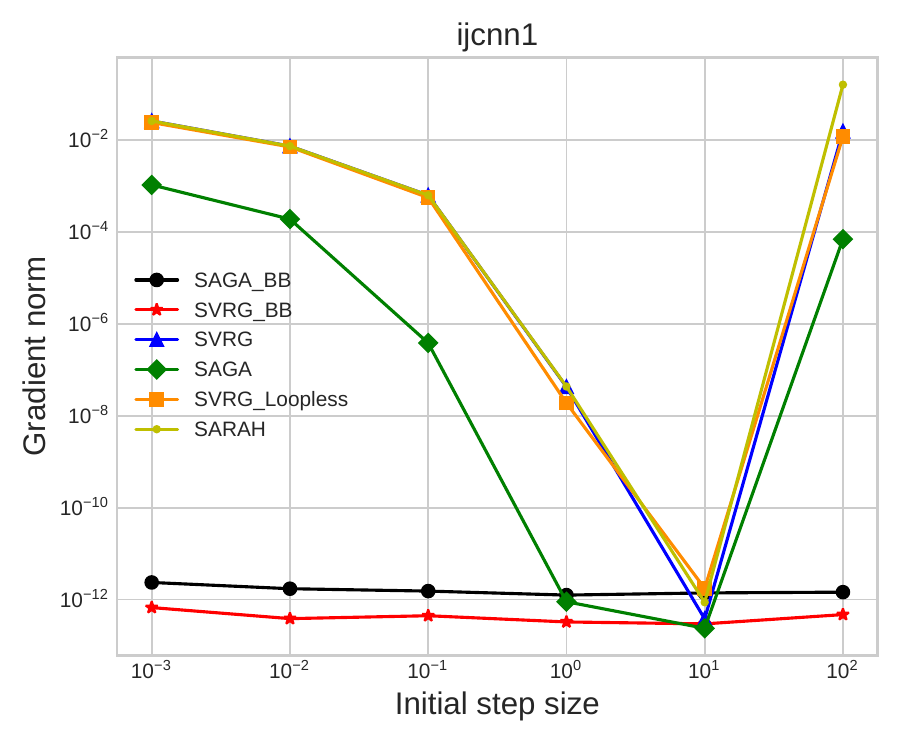}}
    \subfigure[]{
    \includegraphics[width=0.23\textwidth]{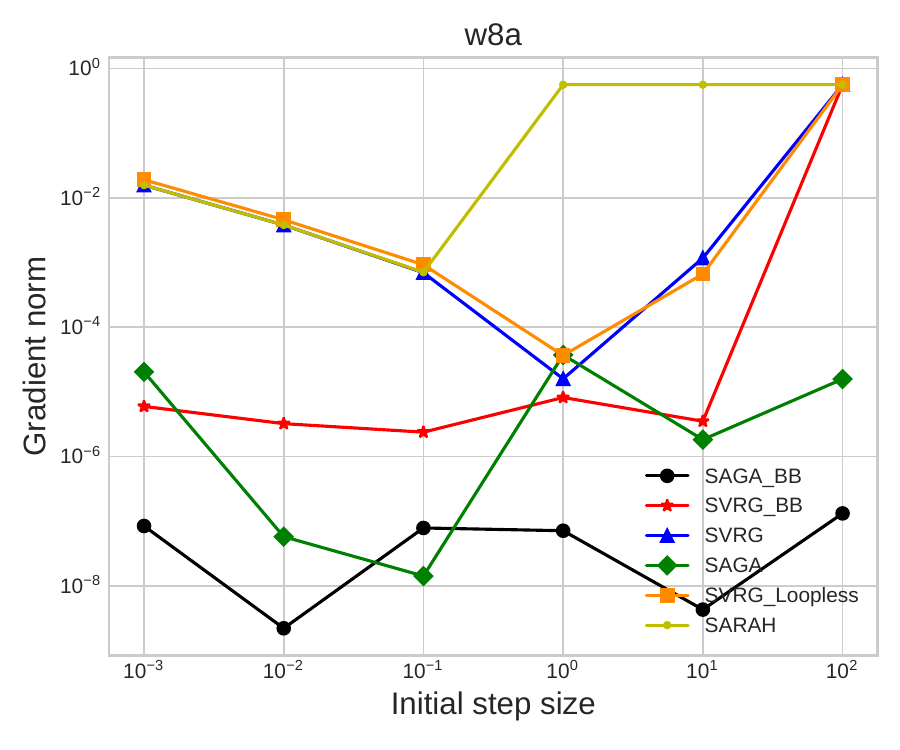}}
    \subfigure[]{
    \includegraphics[width=0.23\textwidth]{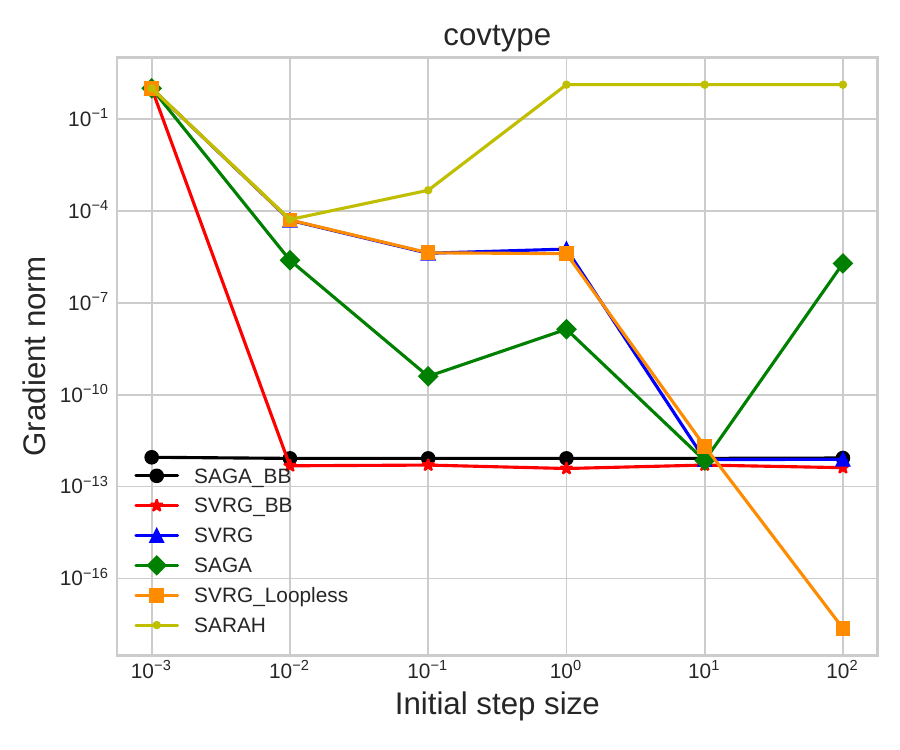}}
    \caption{Performance under different initial step size and fixed batch-size of 16. To make the plot readable, we limit the gradient norm to a maximum value of 1.}
    \label{fig: logistic loss bz16}
\end{figure}

\begin{figure}[ht]
    \centering
    \subfigure[]{
    \includegraphics[width=0.23\textwidth]{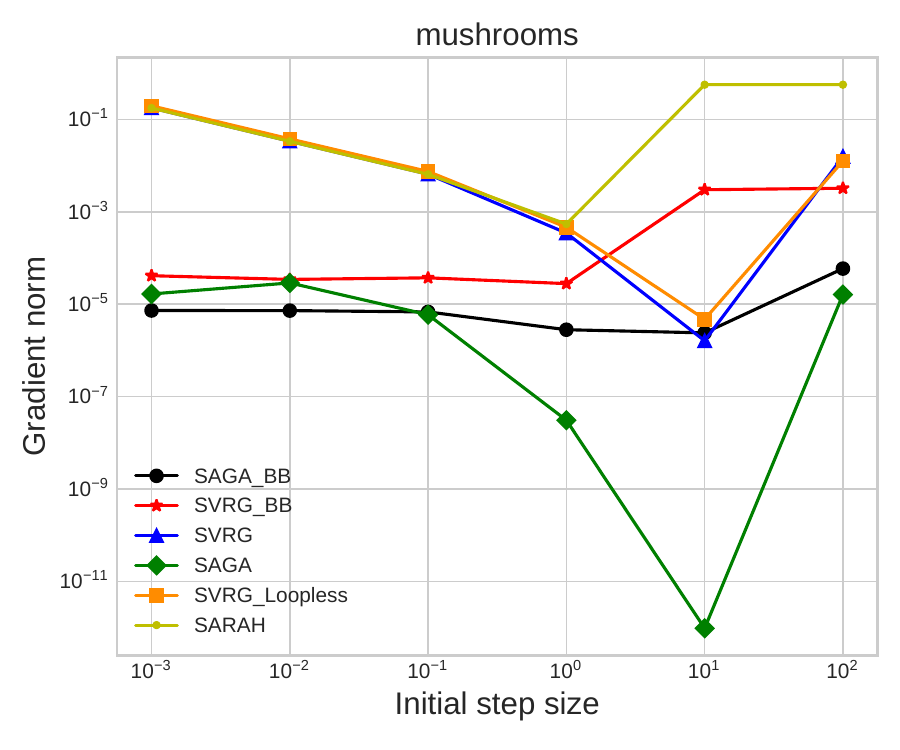}}
    \subfigure[]{
    \includegraphics[width=0.23\textwidth]{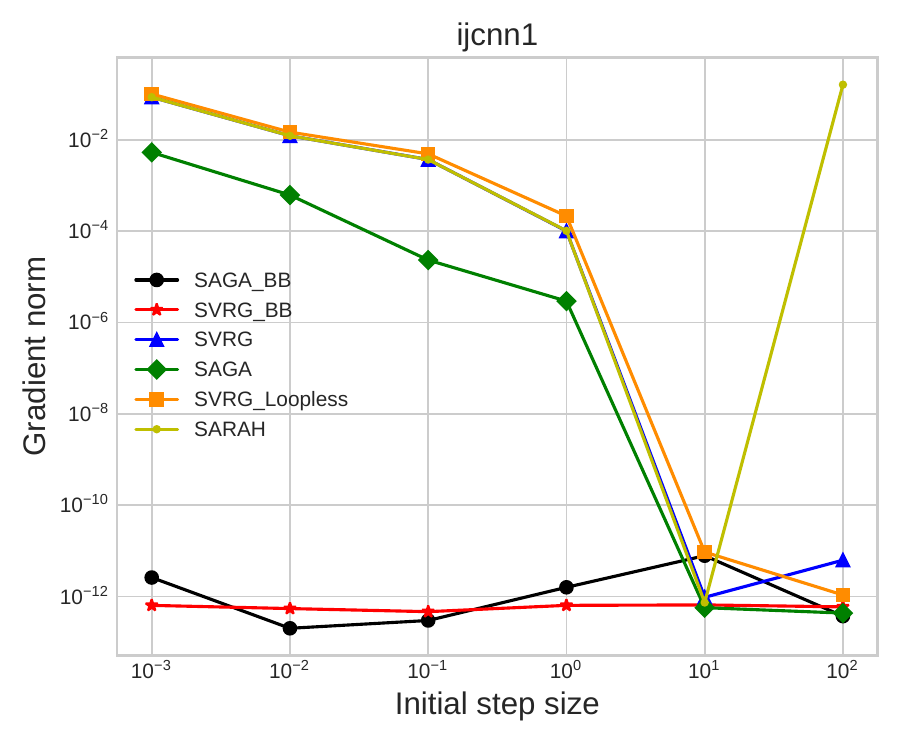}}
    \subfigure[]{
    \includegraphics[width=0.23\textwidth]{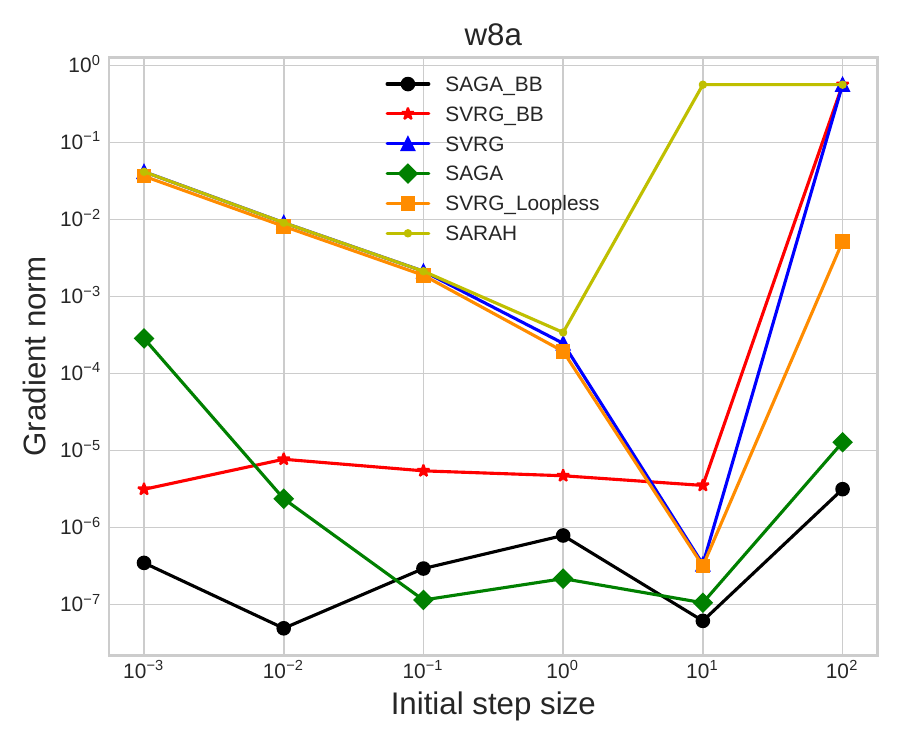}}
    \subfigure[]{
    \includegraphics[width=0.23\textwidth]{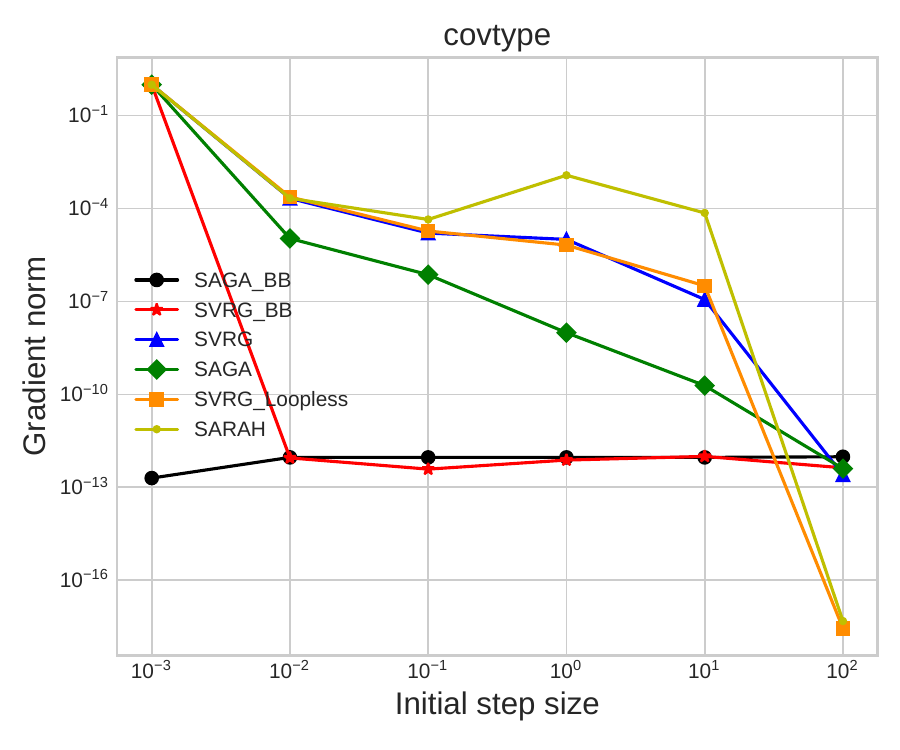}}
    \caption{Performance under different initial step size and fixed batch-size of 64. To make the plot readable, we limit the gradient norm to a maximum value of 1.}
    \label{fig: logistic loss bz64}
\end{figure}

\newpage
\subsection{Additional experiments for Huber loss}\label{sub: Huber loss}
In this subsection, We apply our method and other variance-reduced methods for solving binary $\ell_2$-regularized problems with Huber loss. These figures show that the good performance of SAGA-BB is consistent across losses, batch-sizes and datasets.

\begin{figure}[ht]
    \centering
    \subfigure[]{
    \includegraphics[width=0.23\textwidth]{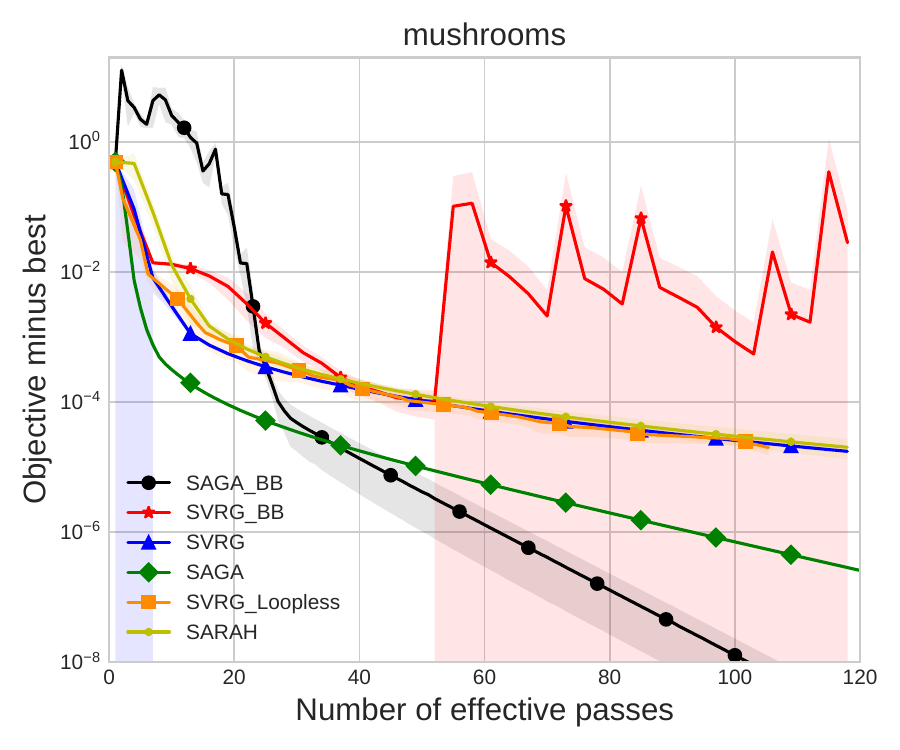}}
    \subfigure[]{
    \includegraphics[width=0.23\textwidth]{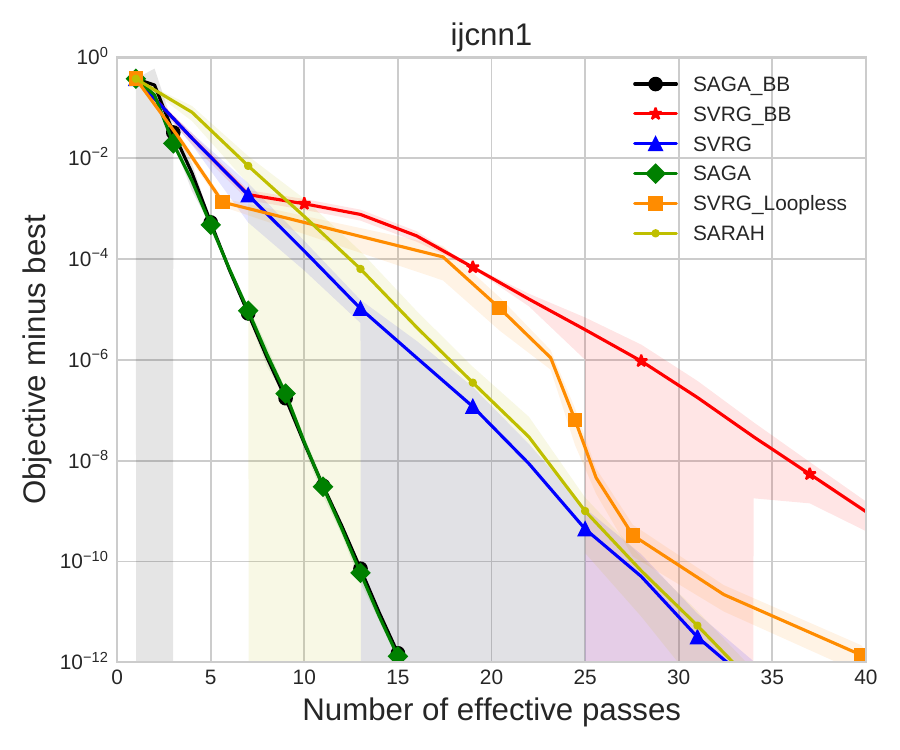}}
    \subfigure[]{
    \includegraphics[width=0.23\textwidth]{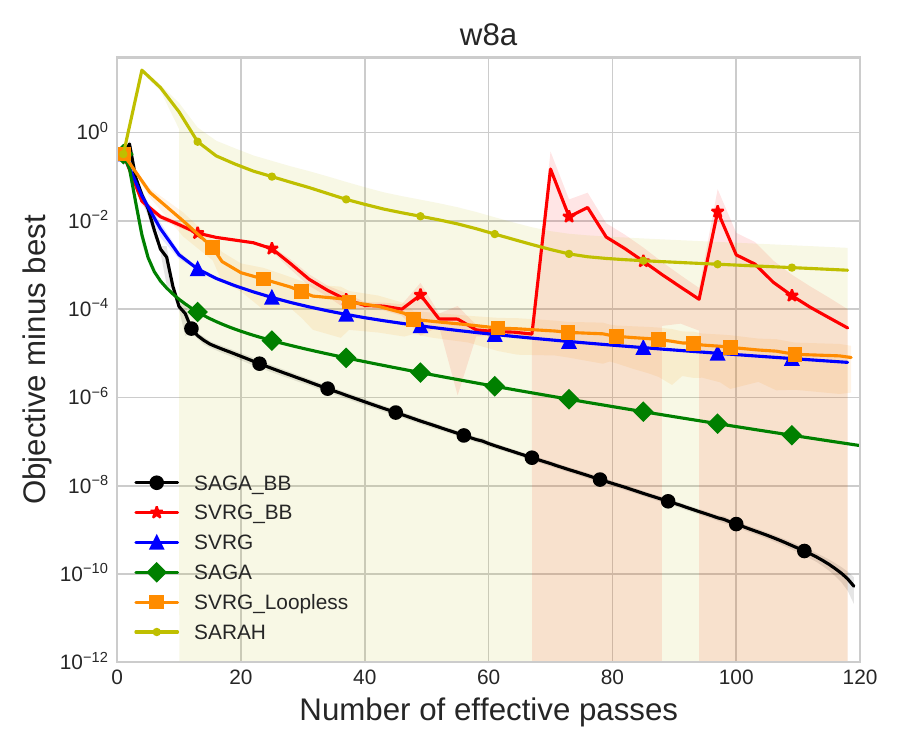}}
    \subfigure[]{
    \includegraphics[width=0.23\textwidth]{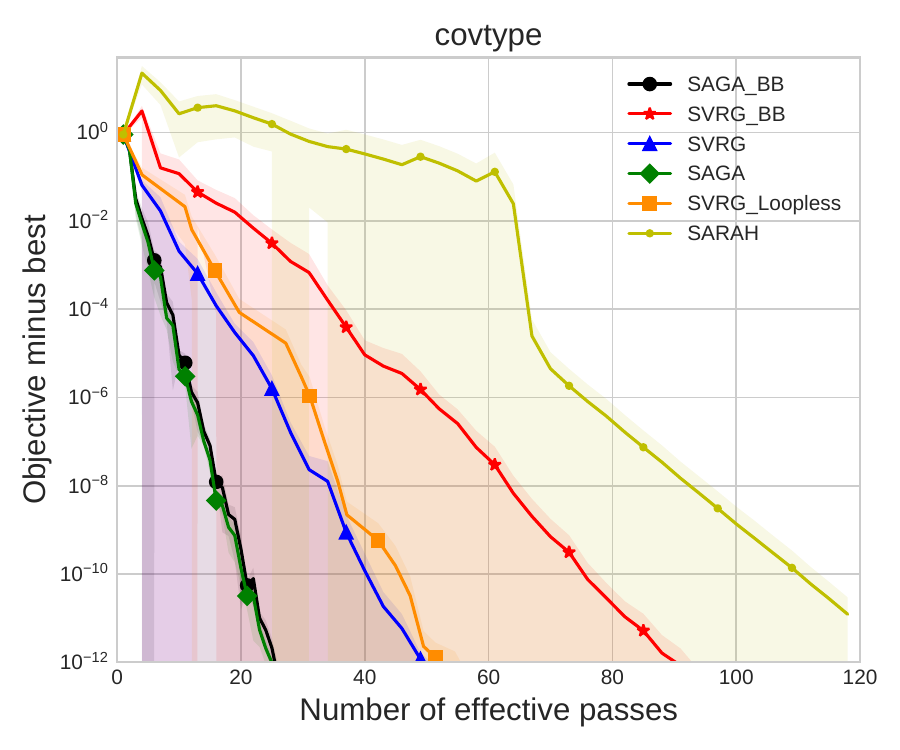}}
    \caption{Convergence on different LibSVM datasets for Huber loss}
    \label{fig: vis huber loss}
\end{figure}

\begin{figure}[ht]
    \centering
    \subfigure[]{
    \includegraphics[width=0.23\textwidth]{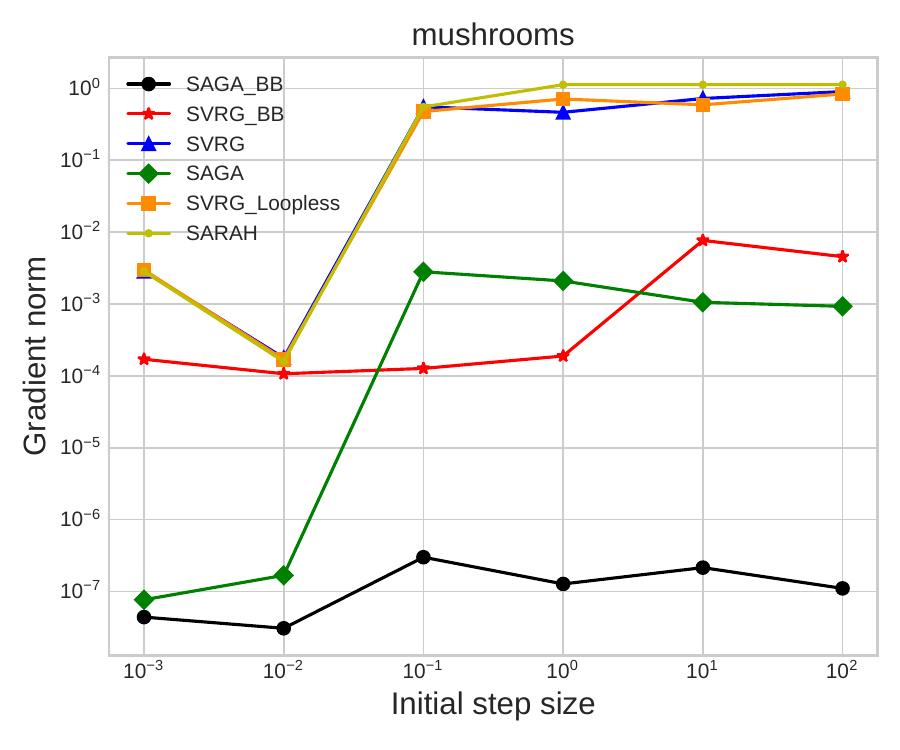}}
    \subfigure[]{
    \includegraphics[width=0.23\textwidth]{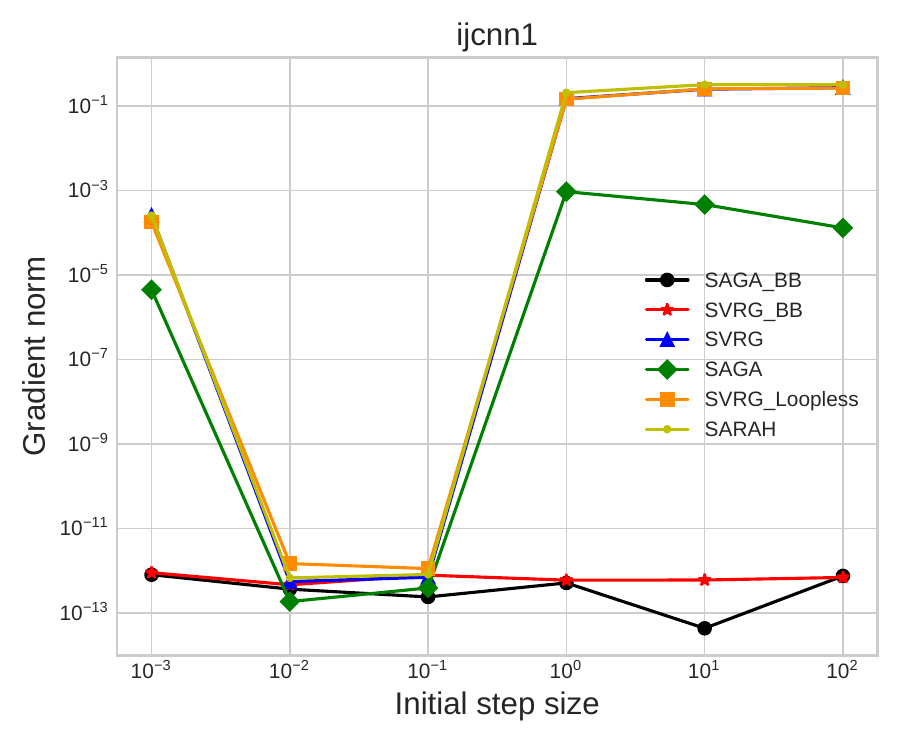}}
    \subfigure[]{
    \includegraphics[width=0.23\textwidth]{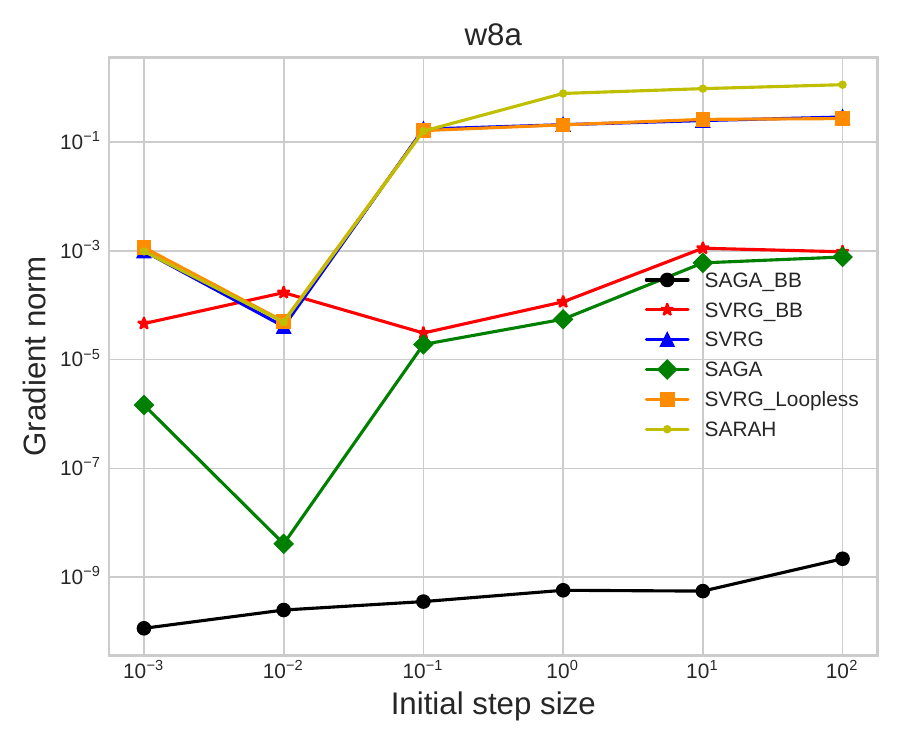}}
    \subfigure[]{
    \includegraphics[width=0.23\textwidth]{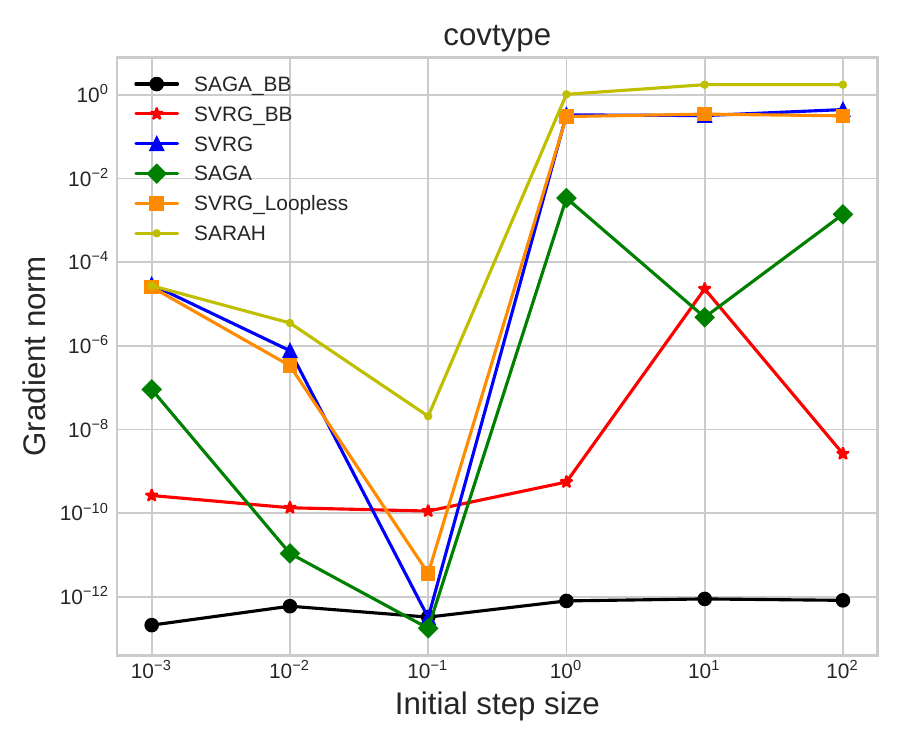}}
    \caption{Performance under different initial step size and fixed batch-size of 1. To make the plot readable, we limit the gradient norm to a maximum value of 1.}
    \label{fig: huber loss bz1}
\end{figure}

\begin{figure}[ht]
    \centering
    \subfigure[]{
    \includegraphics[width=0.23\textwidth]{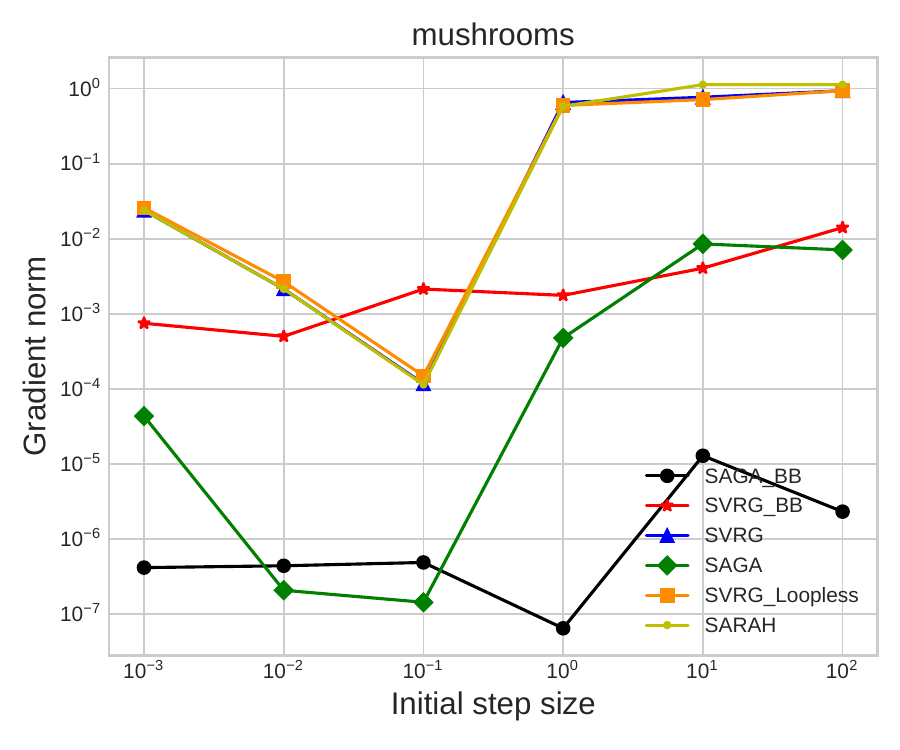}}
    \subfigure[]{
    \includegraphics[width=0.23\textwidth]{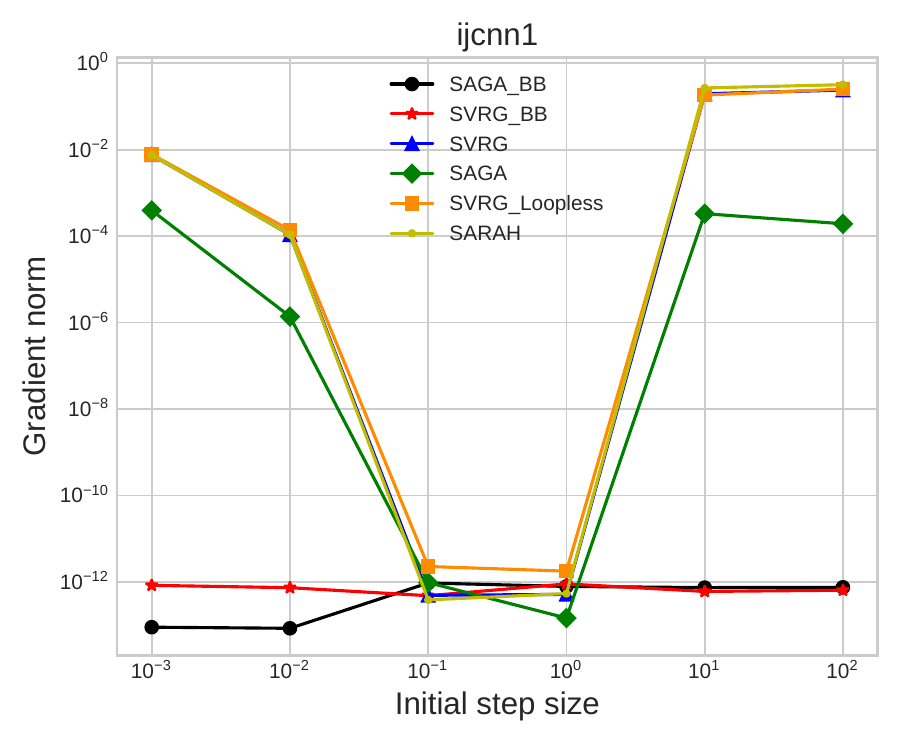}}
    \subfigure[]{
    \includegraphics[width=0.23\textwidth]{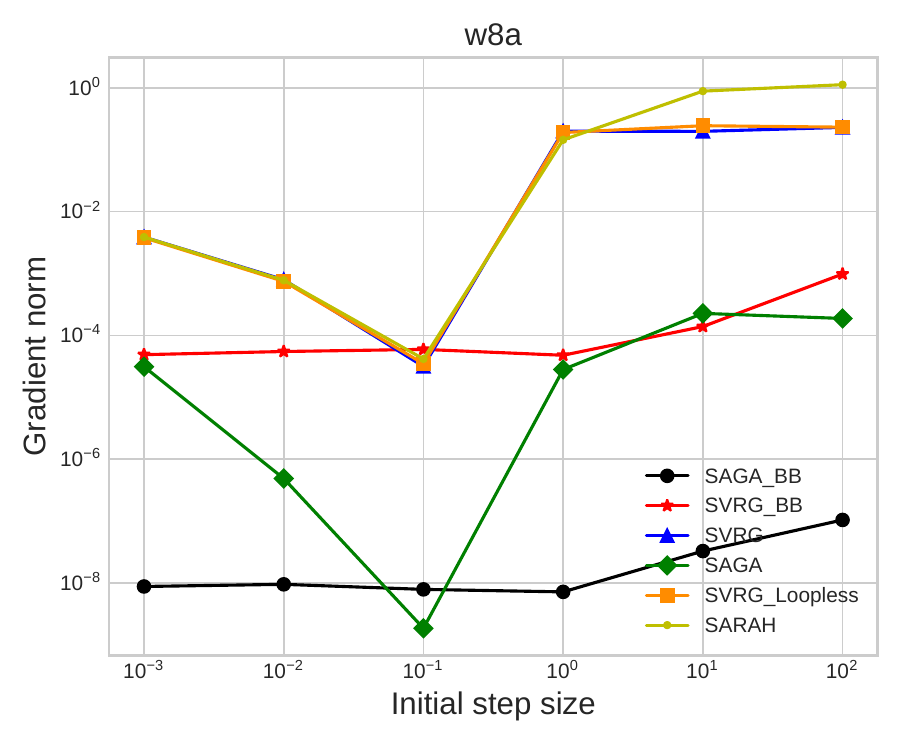}}
    \subfigure[]{
    \includegraphics[width=0.23\textwidth]{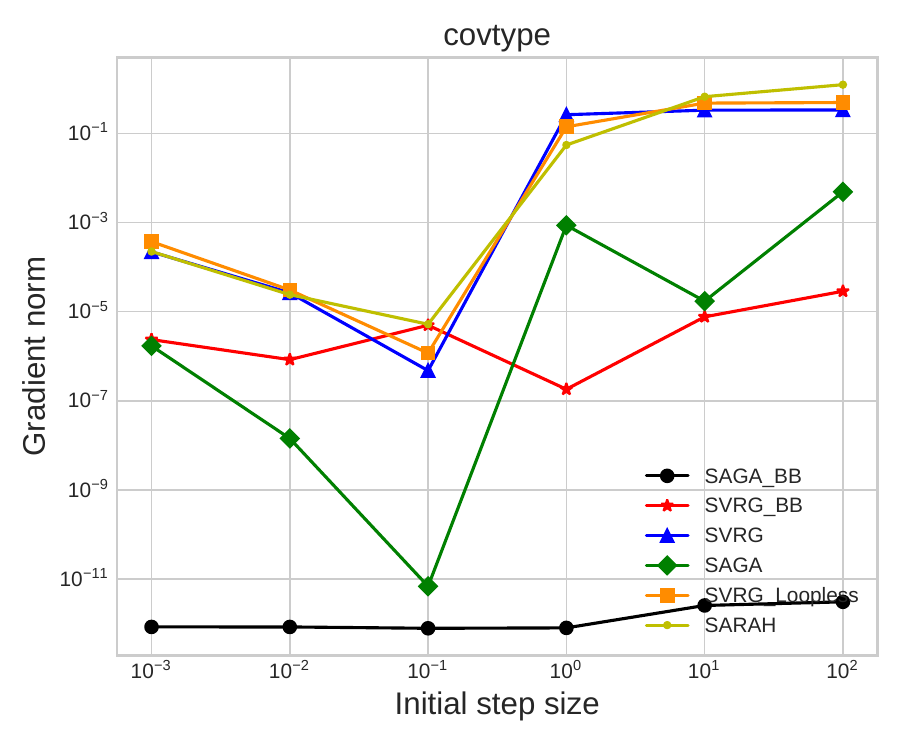}}
    \caption{Performance under different initial step size and fixed batch-size of 8. To make the plot readable, we limit the gradient norm to a maximum value of 1.}
    \label{fig: huber loss bz8}
\end{figure}

\begin{figure}[ht]
    \centering
    \subfigure[]{
    \includegraphics[width=0.23\textwidth]{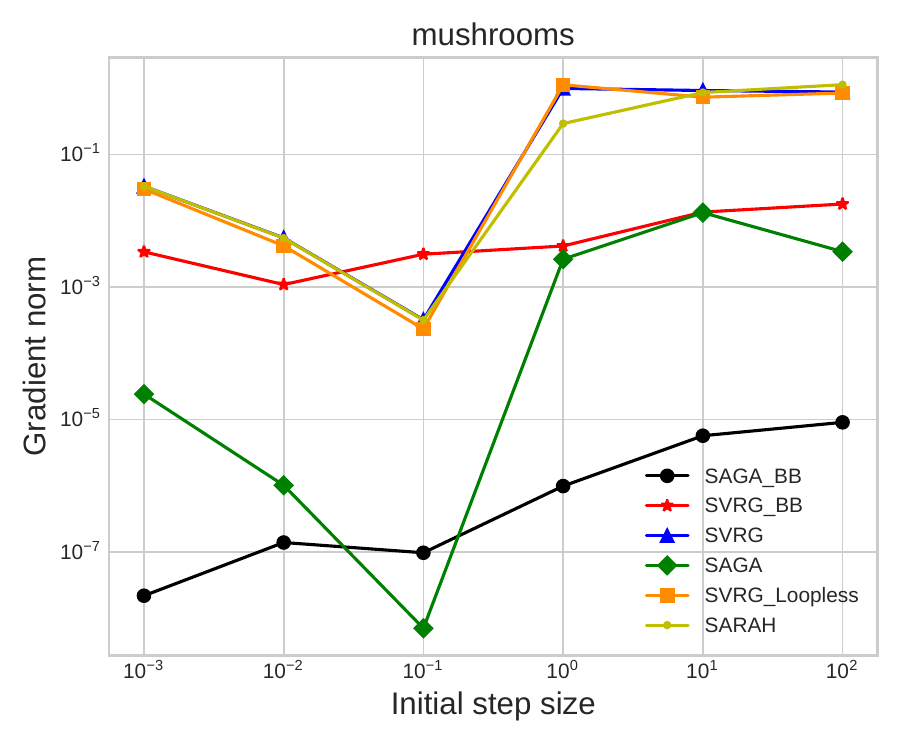}}
    \subfigure[]{
    \includegraphics[width=0.23\textwidth]{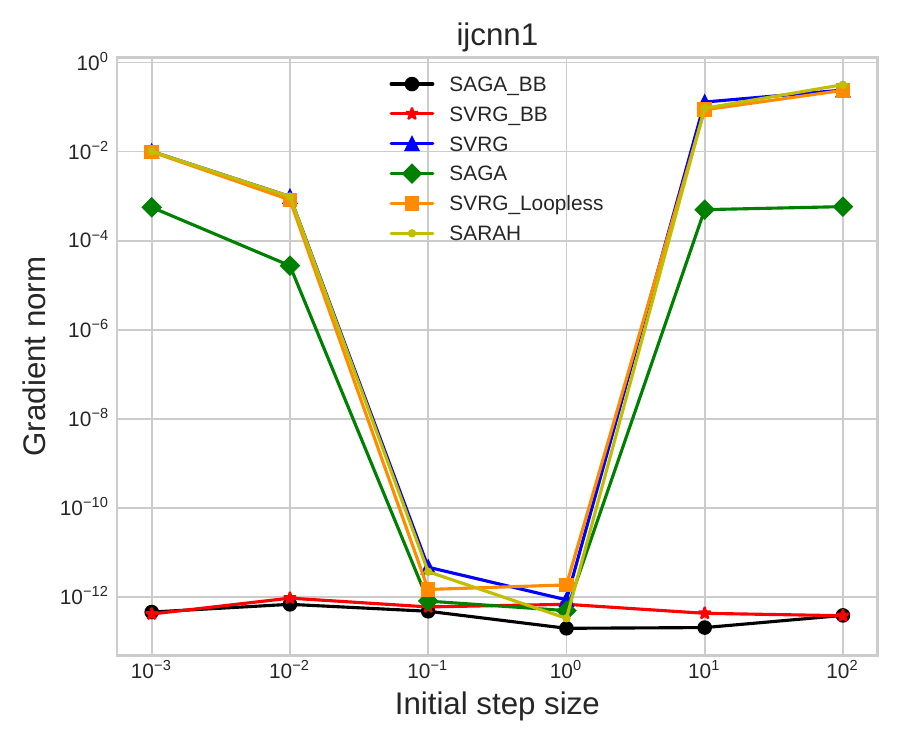}}
    \subfigure[]{
    \includegraphics[width=0.23\textwidth]{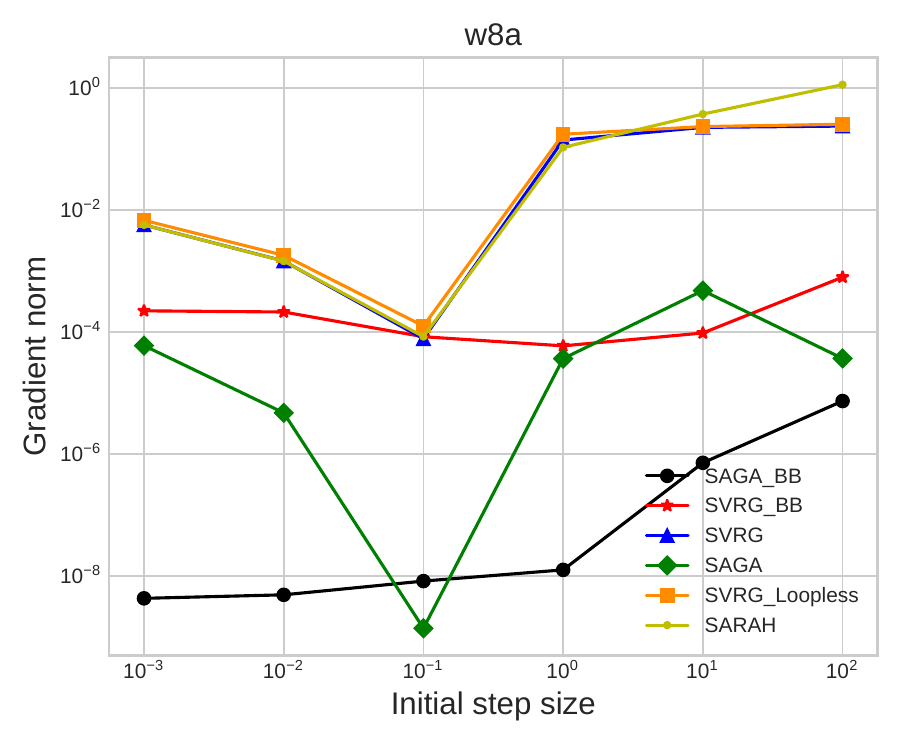}}
    \subfigure[]{
    \includegraphics[width=0.23\textwidth]{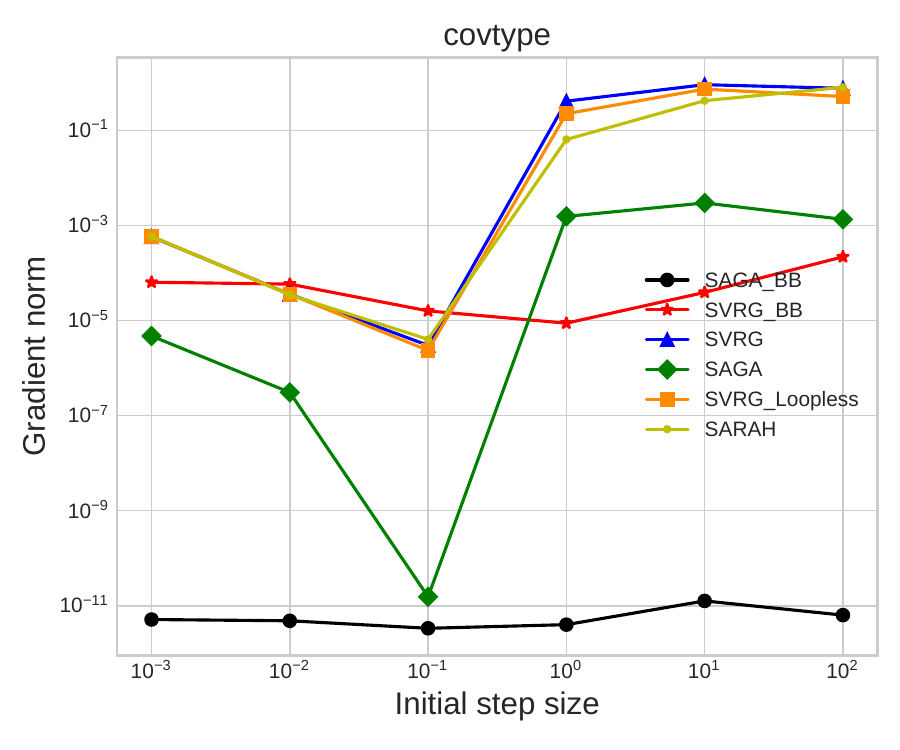}}
    \caption{Performance under different initial step size and fixed batch-size of 16. To make the plot readable, we limit the gradient norm to a maximum value of 1.}
    \label{fig: huber loss bz16}
\end{figure}

\begin{figure}[ht]
    \centering
    \subfigure[]{
    \includegraphics[width=0.23\textwidth]{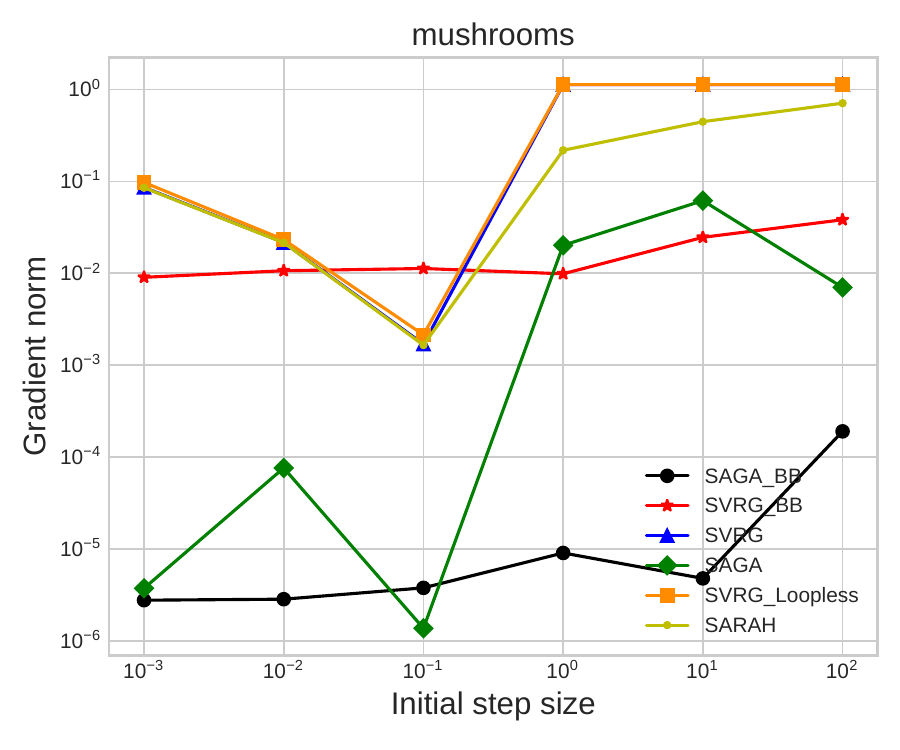}}
    \subfigure[]{
    \includegraphics[width=0.23\textwidth]{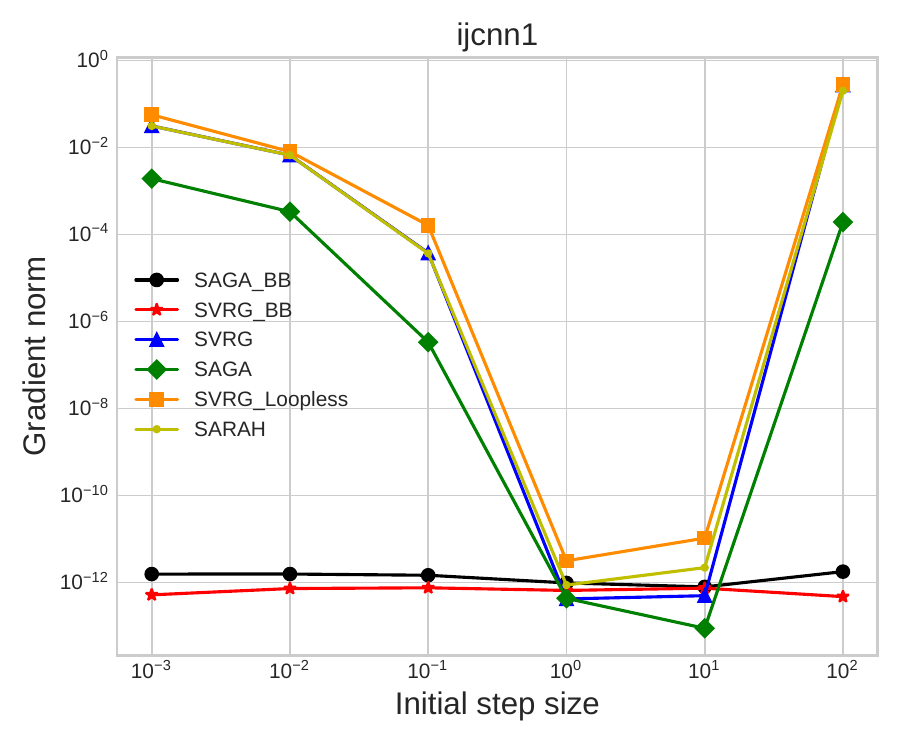}}
    \subfigure[]{
    \includegraphics[width=0.23\textwidth]{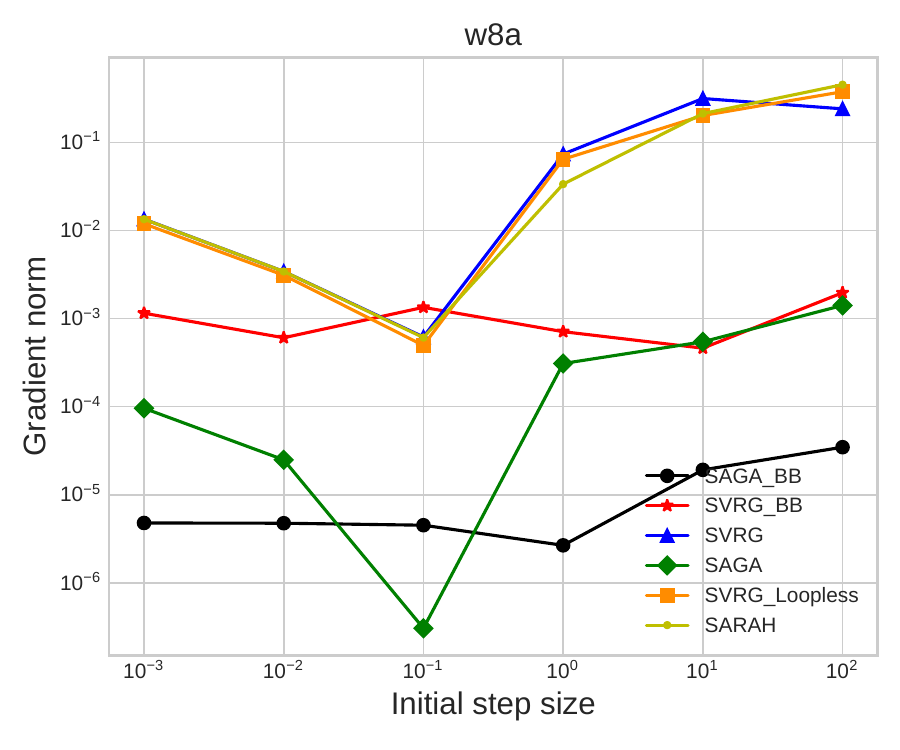}}
    \subfigure[]{
    \includegraphics[width=0.23\textwidth]{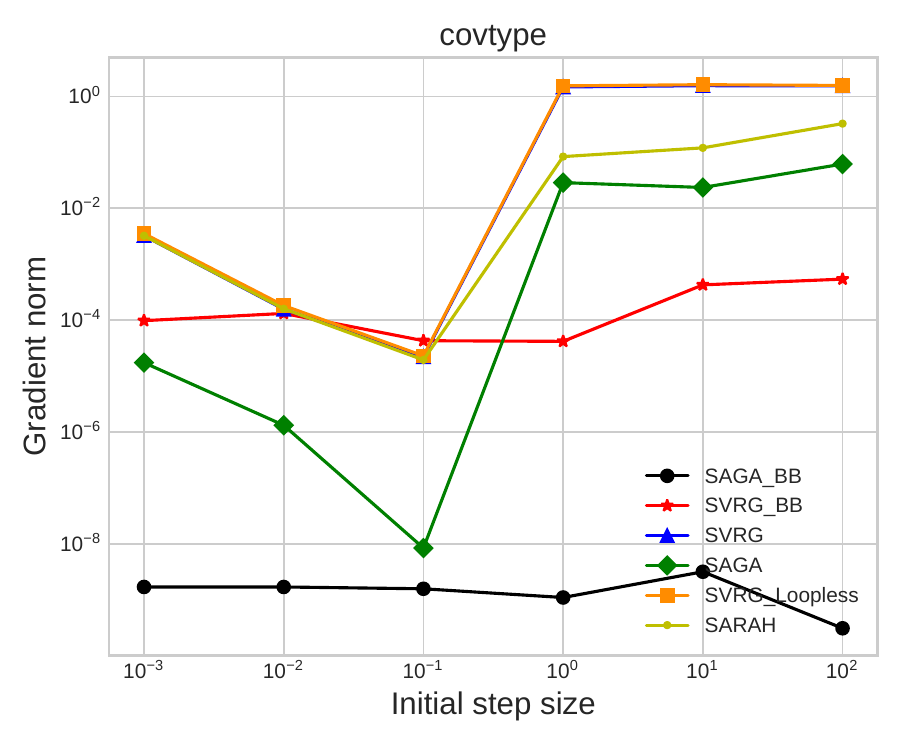}}
    \caption{Performance under different initial step size and fixed batch-size of 64. To make the plot readable, we limit the gradient norm to a maximum value of 1.}
    \label{fig: huber loss bz64}
\end{figure}

\vfill

\end{document}